\documentclass{article}
\usepackage{amsmath,amsrefs,amssymb,amsthm}
\usepackage{calc}
\usepackage[shortlabels]{enumitem}
\usepackage[T1]{fontenc}
\usepackage{forest}
\usepackage{geometry}
\usepackage[colorlinks = true]{hyperref}
\usepackage{cleveref}
\usepackage{mathtools}
\usepackage{parskip}
\usepackage{tikz}
\usepackage{xcolor}

\newtheorem{theorem}{Theorem}[section]

\newtheorem{corollary}[theorem]{Corollary}
\newtheorem{lemma}[theorem]{Lemma}
\newtheorem{proposition}[theorem]{Proposition}

\theoremstyle{definition}
\newtheorem{definition}[theorem]{Definition}
\newtheorem{remark}[theorem]{Remark}
\newtheorem{example}[theorem]{Example}

\newcommand{\N}{\mathbb{N}}

\newcommand{\C}{\mathbb{C}}
\newcommand{\Ht}{\widetilde{H}}

\newcommand{\LD}{\mathsf{LD}}
\newcommand{\inv}{\mathsf{inv}}
\newcommand{\sminv}{\mathsf{sminv}}
\newcommand{\dinv}{\mathsf{dinv}}
\newcommand{\sdinv}{\mathsf{sdinv}}
\newcommand{\area}{\mathsf{area}}
\newcommand{\maj}{\mathsf{maj}}
\newcommand{\pmaj}{\mathsf{pmaj}}
\newcommand{\minimaj}{\mathsf{minimaj}}
\newcommand{\height}{\mathsf{ht}}
\renewcommand{\th}{^\mathsf{th}}
\newcommand{\start}{\mathsf{start}}
\newcommand{\st}{\mathsf{st}}
\newcommand{\SW}{\mathsf{SW}}
\newcommand{\SF}{\mathsf{SF}}
\newcommand{\OP}{\mathsf{OP}}
\newcommand{\Split}{\mathsf{Split}}

\newcommand{\skel}{\mathsf{skel}}

\newcommand{\ul}[1]{\underline{#1}}

\DeclareMathOperator{\grFrob}{grFrob}
\DeclareMathOperator{\Frob}{Frob}

\newcommand{\qbinom}[3][q]{\genfrac{[}{]}{0pt}{}{#2}{#3}_{#1}}

\title{Smirnov words and the Delta Conjectures}
\author{
    Alessandro Iraci\footnote{Università di Pisa, \href{mailto:alessandro.iraci@unipi.it}{alessandro.iraci@unipi.it}}, 
    Philippe Nadeau\footnote{CNRS, Université Lyon 1, \href{mailto:nadeau@math.univ-lyon1.fr}{nadeau@math.univ-lyon1.fr}. P.N. is partially supported by the project ANR19-CE48-011-01 (COMBIN\'E)}, 
    Anna Vanden Wyngaerd\footnote{Université Libre de Bruxelles, \href{emailto:anna.vanden.wyngaerd@ulb.be}{anna.vanden.wyngaerd@ulb.be}}}

\date{\today}

\begin{document}
\maketitle
\begin{abstract}
    We provide a combinatorial interpretation of the symmetric function $\left.\Theta_{e_k}\Theta_{e_l}\nabla e_{n-k-l}\right|_{t=0}$ in terms of segmented Smirnov words. The motivation for this work is the study of a diagonal coinvariant ring with one set of commuting and two sets of anti-commuting variables, whose Frobenius characteristic is conjectured to be the symmetric function in question. Furthermore, this function is related to the Delta conjectures. Our work is a step towards a unified formulation of the two versions, as we prove a unified Delta theorem at $t=0$.
\end{abstract}

\section*{Introduction}

In the 1990s, Garsia and Haiman introduced the ring of \emph{diagonal coinvariants} $DR_n$. The study of the structure of this $\mathfrak S_n$-module and its generalizations has been an important research topic in algebra and combinatorics ever since. There is in fact a rich interplay at work here.  Various conjectural enumerative results on the ring (graded dimensions, multiplicities, etc.) indicate some beautiful underlying combinatorial structures. Understanding the structures leads to information on the ring, while generalizing them raises new questions on the algebraic side, leading in particular to the introduction of various new coinvariant rings. Then, some new combinatorics may emerge from those, and so on. This work will be motivated by a particular coinvariant ring in this realm, for which we unearth and explore the conjectural combinatorics in detail. This in turn gives in particular new insight on the algebraic side. Before we state our results, we need to briefly recall the story of $DR_n$ and its extensions.

The ring $DR_n$ is defined as follows: consider the space $\C[\mathbf{x}_n,\mathbf{y}_n]\coloneqq \C[x_1,\dots,x_n,y_1,\dots,y_n]$ and define an $\mathfrak S_n$-action as \[\sigma\cdot f(x_1,\dots,x_n, y_1,\dots, y_n) \coloneqq f(x_{\sigma(1)},\dots,x_{\sigma(n)}, y_{\sigma(1)},\dots, y_{\sigma(n)})\] for all $f\in \C[\mathbf{x}_n,\mathbf{y}_n]$ and $\sigma \in \mathfrak S_n$. Let $I(\mathbf{x}_n, \mathbf{y}_n)$ be the ideal generated by the $\mathfrak S_n$-invariants with vanishing constant term. Then the ring of diagonal coinvariants is defined as \[DR_n \coloneqq \C[\mathbf{x}_n,\mathbf{y}_n]/I(\mathbf{x}_n,\mathbf{y}_n).\]

This space has a natural bi-grading: let $DR_n^{(i,j)}$ be the component of $DR_n$ with homogeneous $\mathbf{x}$-degree $i$ and homogeneous $\mathbf{y}$-degree $j$. This grading is preserved by the $\mathfrak S_n$-action. Garsia and Haiman conjectured, and Haiman later proved \cite{Haiman2002HilbertScheme}, a formula for the graded Frobenius characteristic of the diagonal harmonics: \[ \grFrob(DR_n; q,t) \coloneqq \sum_{i,j\in\N}q^it^j\Frob(DR_n^{(i,j)})= \nabla e_n,\] where $e_n$ is the $n\th$ elementary symmetric function and $\nabla$ is the operator depending on $q,t$ introduced in \cite{BergeronGarsiaHaimanTesler1999IdentitiesPositivityConjectures}.

In \cite{HaglundHaimanLoehrRemmelUlyanov2005ShuffleConjecture}, the authors gave a combinatorial formula for this graded Frobenius character $\nabla e_n$, in terms of \emph{labelled Dyck paths}, called the \emph{shuffle conjecture}, which is now a theorem, by Carlsson and Mellit \cite{CarlssonMellit2018ShuffleConjecture}. The \emph{Delta conjecture} is a pair of combinatorial formulas for the symmetric function $\Delta'_{e_{n-k-1}}e_n$ in terms of \emph{decorated} labelled Dyck paths, stated in \cite{HaglundRemmelWilson2018DeltaConjecture}, that reduces to the shuffle theorem when $k=0$ (see \eqref{eq:delta-rise} and \eqref{eq:delta-valley}). Here $\Delta'_{e_j}$ is another operator on symmetric functions.

This extension of the combinatorial setting led Zabrocki and, later on, D'Adderio, Iraci and Vanden Wyngaerd to introduce extensions of $DR_n$ \cites{Zabrocki2019ModuleDeltaConjecture,DAdderioIraciVandenWyngaerd2021ThetaOperators}. Consider the ring \[\C[\mathbf{x}_n,\mathbf{y}_n, \boldsymbol{\theta}_n, \boldsymbol{\xi}_n] \coloneqq \C[x_1,\dots,x_n,y_1,\dots,y_n,\theta_1,\dots,\theta_n,\xi_1,\dots,\xi_n],\] where the $\mathbf{x}_n, \mathbf{y}_n$ are the usual commuting variables, and the $\boldsymbol{\theta}_n$, $\boldsymbol{\xi}_n$ are anti-commuting:
\begin{center}
    $\theta_i \theta_j = -\theta_j \theta_i$ and $\xi_i \xi_j = -\xi_j \xi_i$ for all $1 \leq i,j \leq n$.
\end{center}

Borrowing terminology from supersymmetry, the commuting and anti-commuting variables are sometimes referred to as \emph{bosonic} and \emph{fermionic}, respectively. Again, consider the $\mathfrak S_n$-action that permutes all the variables simultaneously:
\begin{align*}
    \sigma\cdot & f(x_1,\dots,x_n, y_1,\dots, y_n,\theta_1,\dots, \theta_n,\xi_1,\dots,\xi_n) \\ &\coloneqq f(x_{\sigma(1)},\dots,x_{\sigma(n)}, y_{\sigma(1)},\dots, y_{\sigma(n)},\theta_{\sigma(1)},\dots, \theta_{\sigma(n)},\xi_{\sigma(1)},\dots,\xi_{\sigma(n)})
\end{align*} for all $f\in \C[\mathbf{x}_n,\mathbf{y}_n, \boldsymbol{\theta}_n, \boldsymbol{\xi}_n]$ and $\sigma \in \mathfrak S_n$. If $I(\mathbf{x}_n,\mathbf{y}_n, \boldsymbol{\theta}_n, \boldsymbol{\xi}_n)$ again denotes the ideal generated by the invariants of this action without constant term, define
\[TDR_n \coloneqq \C[\mathbf{x}_n,\mathbf{y}_n, \boldsymbol{\theta}_n, \boldsymbol{\xi}_n]/I(\mathbf{x}_n,\mathbf{y}_n, \boldsymbol{\theta}_n, \boldsymbol{\xi}_n),\] where the `T' stands for `Twice' as many generators. This ring is naturally quadruply graded: let $TDR_n^{(i,j,k,l)}$ denote the component of $TDR_n$ of homogeneous $(i,j,k,l)$-degrees in the four different sets of variables $(\mathbf{x}_n,\mathbf{y}_n, \boldsymbol{\theta}_n, \boldsymbol{\xi}_n)$.

We sometimes call this the $(2,2)$ case, referring to the 2 sets of commuting and 2 sets of anti-commuting variables. All these modules are actually special cases of a general family of bosonic-fermionic diagonal coinvariant modules, defined by the same construction with $a$ sets of bosonic variables and $b$ sets of fermionic variables. These modules, which are actually $\mathfrak S_n \times \mathsf{GL}_a \times \mathsf{GL}_b$ modules with the action of the general linear group on the matrices of variables, have been studied extensively and are tied to some very nice combinatorial and algebraic results. See for example \cite{Bergeron2020} for a survey on the topic and the statement of the so-called combinatorial supersymmetry.

In \cite{Zabrocki2019ModuleDeltaConjecture} Zabrocki conjectured that, in the $(2,1)$-case of two sets of bosonic and one set of fermionic variables
\begin{equation}
    \label{eq:Delta-module}\grFrob \left( \bigoplus_{i,j\in\N} TDR_n^{(i,j,k,0)} \right) \coloneqq \sum_{i,j\in\N} q^{i} t^{j} \Frob(TDR_n^{(i,j,k,0)}) \overset{?}{=} \Delta'_{e_{n-k-1}} e_n,
\end{equation}
where we recognize on the right-hand size the symmetric function of the Delta conjecture. In \cite{DAdderioIraciVandenWyngaerd2021ThetaOperators}, D'Adderio, Iraci and Vanden Wyngaerd introduced the symmetric function operator $\Theta_f$, for any symmetric function $f$, and showed that $\Delta'_{n-k-1} e_n = \Theta_{e_k} \nabla e_{n-k}$. This new operator permitted them to extend Zabrocki's conjecture to the $(2,2)$-case:
\begin{equation}
    \label{eq:thetatheta-module}
    \grFrob \left( \bigoplus_{i,j\in\N}TDR_n^{(i,j,k,l)} \right) \coloneqq \sum_{i,j\in\N} q^{i} t^{j} \Frob(TDR_n^{(i,j,k,l)}) \overset{?}{=} \Theta_{e_l} \Theta_{e_k} \nabla e_{n-k-l}.
\end{equation}
This general conjecture is still open, but many special cases have been studied over the years. Other than the $(1,0)$- and $(0,1)$- cases, which are classical, there are several recent results. The $(2,0)$-case, or the classical diagonal coinvariant case was proven by Haiman to be a consequence of his famous $n!$-theorem \cite{Haiman2002HilbertScheme} as recalled above. The $(0,2)$-case, or \emph{fermionic Theta} case, was proved by Iraci, Rhoades, and Romero in \cite{IraciRhoadesRomero2023FermionicTheta}. The $(1,1)$-case, or the \emph{superspace coinvariant ring}, is still open; Rhoades and Wilson in \cite{RhoadesWilson2023Superspace} showed that its Hilbert series agrees with the expected formula.\medskip

In this paper, we will turn our interest to the combinatorics of $(1,2)$-case. In other words, we are led through the conjecture~\eqref{eq:thetatheta-module} to study the symmetric function
\begin{equation}
    \label{eq:theSF}
    \left.\Theta_{e_k}\Theta_{e_l}\nabla e_{n-k-l}\right|_{t=0}.
\end{equation}

Our main result, \Cref{thm:main}, is to give a combinatorial interpretation of this formula that is completely new, in terms of \emph{segmented Smirnov words} and a pleasant $q$-statistic. We will also give a variant of this result which deserves to be dubbed ``unified Delta theorem at $t=0$'', based on a $q$-statistic that is more involved. We also rewrite the combinatorial expansion using fundamental quasisymmetric functions, and apply this to extract the coefficient of $s_{1^n}$ in \eqref{eq:theSF}, which is conjecturally the Hilbert series of the sign-isotypic component of the $(1,2)$-case.
\smallskip

Let us give more details about these results. A Smirnov word is a word in the alphabet of nonnegative integers such that no two adjacent letters are the same; see~\cite{ShareshianWachs2016Chromatic} for instance. A segmented Smirnov word is a concatenation of Smirnov words with prescribed lengths (see \Cref{def:smirnov}). Then \Cref{thm:main} is a monomial expansion for the symmetric function \eqref{eq:theSF} as a generating function of segmented Smirnov words, where the exponent of $q$ is given by a new \emph{sminversion} statistic on these words (see \Cref{def:smirnov_inv}). In our formula, the values of $k$ and $l$ in \eqref{eq:theSF} correspond to the number of ascents and descents in the Smirnov word. Our proof relies on an algebraic recursion (\Cref{thm:recursion}) which we show to be a consequence of a symmetric function identity in \cite{DAdderioRomero2023ThetaIdentities}*{Theorem~8.2}. We then show in \Cref{sec:combinatorial_expansion} that the combinatorial side satisfies the same recursion, from which the theorem follows.

Now the full symmetric function on the right hand side of~\eqref{eq:thetatheta-module} is of particular combinatorial interest, as it may provide a unified version of the two different Delta conjectures:
\begin{align}
    \Delta'_{e_{n-k-1}} e_n = \Theta_{e_k} \nabla e_{n-k}
     & = \sum_{D \in \LD(n)^{\ast k}} q^{\dinv(D)} t^{\area(D)} x^D \label{eq:delta-rise}       \\
     & = \sum_{D \in \LD(n)^{\bullet k}} q^{\dinv(D)} t^{\area(D)} x^D \label{eq:delta-valley}.
\end{align}
The sets $\LD(n)^{\ast k}$ and $\LD(n)^{\bullet k}$ denote labelled Dyck paths of size $n$ with $k$ decorations on \emph{rises} or \emph{valleys}, respectively; and the statistics $\dinv$ and $\area$ depend on the decorations (see \Cref{sec:link-delta} for precise definitions).
So \eqref{eq:delta-rise} is referred to as the \emph{rise version} and \eqref{eq:delta-valley} as the \emph{valley version} of the Delta conjecture. The rise version was recently proved in \cite{DAdderioMellit2022CompositionalDelta}.

In \cite{DAdderioIraciVandenWyngaerd2021ThetaOperators}, the first and third named authors of this paper, together with Michele D'Adderio, conjectured a partial formula for a possible \emph{unified Delta conjecture}, for which they have significant computational evidence:
\begin{equation}
    \left. \Theta_{e_k} \Theta_{e_l} \nabla e_{n-k-l} \right|_{q=1} = \sum_{D \in \LD(n)^{\ast k, \bullet l}} t^{\area(D)} x^D \label{eq:unified-delta-q=1},
\end{equation}
where $\LD(n)^{\ast k, \bullet l}$ is the set of labelled Dyck paths of size $n$ with $k$ decorations on rises and $l$ decorations on valleys. In order to get fully unified Delta conjecture, one has to find a statistic $\mathsf{qstat} \colon \LD(n)^{\ast k, \bullet l} \rightarrow \N$ so that
\begin{equation}
    \Theta_{e_k} \Theta_{e_l} \nabla e_{n-k-l} = \sum_{D \in \LD(n)^{\ast k, \bullet l}} q^{\mathsf{qstat}(D)} t^{\area(D)} x^D \label{eq:unified-delta}
\end{equation}
and such that when $k=0$ or $l = 0$, the formula would reduce to \eqref{eq:delta-rise} or \eqref{eq:delta-valley}, respectively.

We give a partial answer to this question in~\Cref{sec:link-delta}, namely we provide such a $q$-statistic at $t=0$: this is~\Cref{thm:unified_Delta_t_equals_0}, our unified Delta theorem at $t=0$. We describe an explicit bijection $\phi$ between segmented Smirnov words and elements of $\LD(n)^{\ast k, \bullet l}$ (doubly decorated labelled Dyck paths) with $\area$ equal to $0$; we then introduce a variant of the $q$-statistic $\sminv$ on segmented Smirnov words, which we call $\sdinv$, so that through the bijection $\phi$ we recover the known $\dinv$ statistic on decorated Dyck paths when $k=0$ or $l=0$, solving \eqref{eq:unified-delta} when $t=0$.

In the last sections, we return to the expansion in~\Cref{thm:main}. In~\Cref{sec:fundamental-quasisym}, we first explicit the fundamental quasisymmetric function expansion of \eqref{eq:theSF} in terms of our combinatorics (\Cref{prop:fundamental}). This is then used to extract the coefficient of $s_{1^n}$ in~\eqref{eq:theSF}, which turns out to have a nice product formula, cf \Cref{prop:catalan-sw}. We finally discuss in \Cref{sec:topdegree} the special case $n=k+l+1$ --corresponding to an expansion in Smirnov words-- that links our work with other areas of combinatorics: the chromatic quasisymmetric function of the path graph (\Cref{sec:chromatic}), parallelogram polyominoes (\Cref{sec:polyominoes}), and noncrossing partitions (\Cref{sec:topdegree}).

\section{Combinatorial definitions}
\label{sec:combinatorial}

In this work $\mathbb{Z}_+$ is the set of positive integers, and we will fix $n\in\mathbb{Z}_+$. We write $\mu \vDash_0 n$ if $\mu$ is a \emph{weak composition} of $n$, that is $\mu=(\mu_1,\mu_2,\ldots)$ where the $\mu_i$ are nonnegative integers that sum to $n$. A composition $\alpha\vDash n$ is a finite sequence $\alpha=(\alpha_1,\ldots,\alpha_s)$ of positive integers that sum to $n$.

\begin{definition}
    \label{def:smirnov}
    A \emph{Smirnov word} of length $n$ is an element $w \in \mathbb{Z}_+^n$ such that $w_i \neq w_{i+1}$ for all $1 \leq i < n$. A \emph{segmented Smirnov word} of shape $\alpha = (\alpha_1, \dots, \alpha_s) \vDash n$ is an element $w \in \mathbb{Z}_+^n$ such that, if we write $w = w^1 w^2 \cdots w^s$ as a concatenation of subwords $w^i$ of length $\alpha_i$, then each $w^i$ is a Smirnov word. We call $w^1, \dots, w^s$ \emph{blocks} of $w$.
\end{definition}

We usually simply denote a segmented Smirnov word by $w$, and omit mention of the shape $\alpha$. In examples, we separate blocks by vertical bars. For instance $ 23|1242|2|31$ is a segmented Smirnov word of length $9$ with shape $(2,4,1,2)$. Note that, if $w$ is a segmented Smirnov word, we allow $w_i = w_{i+1}$ if they belong to different blocks.

Let $\SW(n)$ be the set of segmented Smirnov words of length $n$. Given $\mu \vDash_0 n$, we denote by $\SW(\mu)$ the set of segmented Smirnov words with content $\mu$, meaning that their multiset of letters is $\{i^{\mu_i} \mid i > 0\}$ (i.e.\ they contain $\mu_1$ occurrences of $1$, $\mu_2$ occurrences of $2$, and so on).

We clearly have $\SW(n) = \bigcup_{\mu \vDash_0 n} \SW(\mu)$.

\begin{example}\label{ex:segmented-smirnov-words-21}
    For example, if $\mu = (2,1)$, then $\SW(\mu)$ consists of
    \begin{center}
        \begin{tabular}{llll}
            $1|1|2$ & $1|12$ & $12|1$ & $121$ \\
            $1|2|1$ & $1|21$ & $21|1$         \\
            $2|1|1$ &        &        &
        \end{tabular}
    \end{center}
    where the columns give the words of shape $(1,1,1)$, $(1,2)$, $(2,1)$ and $(3)$, respectively.
\end{example}

It is convenient to introduce the following notation.

\begin{definition}
    Given a Smirnov word $w = w_1 \cdots w_n$, we say that $i$ is
    \begin{itemize}
        \item a \emph{peak} if $w_{i-1} < w_i > w_{i+1}$;
        \item a \emph{valley} if $w_{i-1} > w_i < w_{i+1}$;
        \item a \emph{double rise} if $w_{i-1} < w_i < w_{i+1}$;
        \item a \emph{double fall} if $w_{i-1} > w_i > w_{i+1}$.
    \end{itemize}
\end{definition}

\begin{definition}
    Consider the alphabet $\mathbb{Z}_+ \cup \{\infty\}$ where $\infty$ is larger than all ``finite'' letters. Given a nonempty segmented Smirnov word $w=(w^1, \dots, w^s)$, define the word in the alphabet $\mathbb{Z}_+ \cup\{\infty\}$ as the concatenation \[ a(w) = \infty\;w^1 \infty \; w^2 \infty \; \cdots \; w^{s-1} \infty \; w^s\;\infty. \]

    If $w$ is a segmented Smirnov word, we say that $i$ is a peak, valley, double rise, or double fall of $w$ if it is one in the word $a(w)$.
\end{definition}

We call \emph{segmented permutation} a segmented Smirnov word in $\SW(1^n)$, that is, a segmented Smirnov word whose letters are exactly the numbers from $1$ to $n$. Thus segmented permutations are simply pairs $(\sigma, \alpha)$ with $\sigma \in S_n$ and $\alpha \vDash n$, since the Smirnov condition is automatically satisfied.

\begin{definition}
    Given a Smirnov word $w$, we say that $i$ is an \emph{ascent} of $w$ if $w_{i+1} > w_i$, and a \emph{descent} otherwise. If $w$ is a segmented Smirnov word, we say that $i$ is an ascent (resp.\ descent) of $w$ if it is an ascent (resp.\ descent) of one of its blocks.
\end{definition}

Let us denote by $\SW(n,k,l)$ the set of segmented Smirnov words with exactly $k$ ascents and $l$ descents. Note that such a segmented Smirnov word has exactly $n-k-l$ blocks, as any index $1 \leq i < n$ must be either an ascent, a descent, or the last letter of a block. For $\mu\vDash_0 n$, we also define $\SW(\mu,k,l)$ as the intersection $\SW(\mu)\cap \SW(n,k,l)$.

\begin{definition}
    An index $i \in \{1, \dots, n\}$ is called \emph{initial} (resp.\ \emph{final}) if it corresponds to the first (resp.\ last) position of a block, i.e.\ if it is of the form $i = \alpha_1 + \cdots + \alpha_{m-1} + 1$ (resp.\ $i = \alpha_1 + \cdots + \alpha_m$) for some $m \in \{1, \dots,s\}$.
\end{definition}

We now introduce a statistic on segmented Smirnov words.

\begin{definition}
    \label{def:smirnov_inv}
    For a segmented Smirnov word $w$, we say that $(i,j)$ with $1 \leq i < j \leq n$ is a Smirnov inversion (or \emph{sminversion}) if $w_i > w_j$ and one of the following holds:
    \begin{enumerate}
        \item \label{inv:first_letter} $w_j$ is the first letter of its block;
        \item \label{inv:bigger_entry} $w_{j-1} > w_i$;
        \item \label{inv:equal_then_first_letter} $i \neq j-1$, $w_{j-1} = w_i$, and $w_{j-1}$ is the first letter of its block;
        \item \label{inv:equal_then_bigger_entry} $i \neq j-1$ and $w_{j-2} > w_{j-1} = w_i$.
    \end{enumerate}
    We denote $\sminv(w)$ the number of sminversions of $w$, and \[ \SW_{q}(\mu,k,l) \coloneqq \sum_{w \in \SW(\mu,k,l)} q^{\sminv(w)}. \]
\end{definition}

\begin{remark}
    \label{rem:231}
    Although we gave four cases in this definition in order to be unambiguous, there is a way to think about the sminversions as certain occurrences of $2-31$ pattern. Recall that a usual $2-31$ occurrence is defined for (not necessarily Smirnov) words $a = a_1\cdots a_n$ as a pair $(i,j)$ with $i<j-1$ that satisfies $a_j<a_i<a_{j-1}$.

    It is easy to see that a usual $2-31$ occurrence in $a(w)$ corresponds to a sminversion in $w$ of type \eqref{inv:first_letter} or \eqref{inv:bigger_entry}. If we extend the notion of $2-31$ occurrence to allow $a_i=a_{j-1}$ but only when $a_{j-2}>a_{j-1}$, then these extended occurrences in $a(w)$ are precisely our sminversions in $w$.
\end{remark}

\begin{example}\label{ex:sminv}
    For example, consider the following segmented Smirnov word:
    \begin{center}
        \begin{tabular}{cl}
            the word    & $231|3212|12$                      \\
            the indices & $123\phantom{|}4567\phantom{|}89$.
        \end{tabular}
    \end{center}
    Its sminversions are $(1,3)$, $(1,6)$, $(1,8)$, $(2,5)$, $(2,8)$, $(4,8)$, $(5,8)$ and $(7,8)$; so its $\sminv$ is 8.
\end{example}

Smirnov inversions actually extend inversions on ordered multiset partitions. For $\mu \vDash_0 n$ and $r \in \N$, an ordered multiset partition with content $\mu$ and $r$ blocks is a partition of the multiset $\{i^{\mu_i} \mid i > 0 \}$ into $r$ ordered subsets, called blocks (see for example \cite{Wilson2016Equidistribution}). We denote by $\OP(\mu, r)$ the set of such ordered multiset partitions, and we recall here the definition of a statistic on this set of objects.

\begin{definition}
    \label{def:op_inv}
    For an ordered set partition $\pi = (\pi_1, \dots, \pi_r)$, an inversion is a pair of elements $a > b$ such that $a \in \pi_i$, $b = \min \pi_j$, and $i < j$. We denote by $\inv(\pi)$ the number of inversions of $\pi$.
\end{definition}

From a segmented Smirnov word $w=(w^1,\ldots,w^s)$, we can define $\Pi(w)$ by considering each word $w^i$ only up to permutation of its entries, that is, associating to it the multiset of these entries. For instance $\Pi(231|3212|1212)=\{1,2,3\}|\{1,2,2,3\}|\{1,1,2,2\}$. Note that this is not an ordered multiset partition in general, as blocks of $\Pi(w)$ are not always sets, but we have the following.

\begin{proposition}
    \label{prop:sminv_to_inv}
    The map $\Pi$ restricts to bijections \[ \SW(\mu, k, 0) \leftrightarrow \OP(\mu, n-k), \qquad \SW(\mu, 0, l) \leftrightarrow \OP(\mu, n-l) \] sending $\sminv$ to $\inv$ in both cases.
\end{proposition}

\begin{proof}
    When $k=0$ or $l=0$, the blocks of a segmented Smirnov word are strictly increasing or decreasing, and thus become sets via $\Pi$. Thus, the restrictions of $\Pi$ are well-defined, and it is clear that they are bijections in both cases.

    Let $w \in \SW(\mu, k, 0)$. Since there are no descents, the blocks of $w$ are strictly increasing. Let $(i,j)$ be a sminversion in $w$: then, since the blocks of $w$ are strictly increasing, $j$ is necessarily initial. It is therefore the minimum of its block, and so $w_i > w_j$ is an inversion in the corresponding ordered set partition. Conversely, if $\pi \in \OP(\mu, n-k)$, writing the blocks of $\pi$ in increasing order yields a segmented Smirnov word with no descents, and the inversions of $\pi$ correspond to sminversions for the same reason.

    Let now $w \in \SW(\mu, 0, l)$. Since there are no ascents, the blocks of $w$ are strictly decreasing. Let $(i,j)$ be a sminversion in $w$: then, since the blocks of $w$ are strictly decreasing, $w_j$ is necessarily maximal among the block letters smaller than $w_i$. It follows that $j$ is the unique index in its block forming a sminversion with $i$. Furthermore, in each block to the right of $i$, there is one such index if and only if the final letter of the block, which is also its minimum, is smaller than $w_i$. It follows that each sminversion in $w$ corresponds to exactly one inversion in the corresponding ordered set partition. Conversely, if $\pi \in \OP(\mu, n-l)$, writing the blocks of $\pi$ in decreasing order yields a segmented Smirnov word with no descents, and the inversions of $\pi$ correspond to sminversions for the same reason.
\end{proof}

\section{Algebraic recursion}

For standard notations and conventions for symmetric functions, we refer the reader to \cites{Haglund2008Book,Macdonald1995Book,Stanley1999EnumerativeCombinatoricsV2}. The Theta operators $\Theta_{f}$ are introduced in~\cite{DAdderioIraciVandenWyngaerd2021ThetaOperators}. For example, $\tilde H_\lambda$ denotes the modified Macdonald polynomial of shape $\lambda$. For any symmetric function $f$, the operator $f^\perp$ is such that $\langle f^\perp g, h \rangle = \langle g, fh \rangle$ for all symmetric functions $g,h$; where $\langle\;,\;\rangle$ is the \emph{Hall scalar product} which can be defined by declaring the basis of Schur functions to be orthonormal. The symmetric functions $E_{n,k}$ are defined via the expansion \[e\left[X\frac{1-z}{1-q}\right] = \sum_{k=0}^{n} \frac{(z;q)_k}{(q;q)_k} E_{n,k},\] where we made use of \emph{plethystic notation} and the \emph{$q$-Pochhammer symbol}. Setting $z=q$, we get that $\sum_{k=0}^{n}E_{n,k} = e_n$. The symmetric function $\nabla E_{n,k}$ provides a nice combinatorial refinement of the famous shuffle theorem, in which the $k$ specifies the number of times the Dyck path returns to the diagonal.

Our main algebraic recursion is the specialization $t=0$ of \cite{DAdderioRomero2023ThetaIdentities}*{Theorem~8.2}. We recall here the statement.

\begin{theorem}[{\cite{DAdderioRomero2023ThetaIdentities}*{Theorem~8.2}}]
    \label{thm:full_recursion}
    \begin{align*}
        h_j^\perp & \Theta_{e_m} \Theta_{e_l} \Ht_{(k)}                                                                                                                                                  \\
                  & = \sum_{r=0}^j \qbinom{k}{r} \sum_{a=0}^k \sum_{b=1}^{j-r+a} q^{\binom{k-r-a}{2}} \qbinom{b-1}{a} \qbinom{b+r-a-1}{k-a-1} \Theta_{e_{m-j+r}} \Theta_{e_{l+k-j+a}} \nabla E_{j-r+a,b} \\
                  & + \sum_{r=0}^j \qbinom{k}{r} \sum_{a=0}^k \sum_{b=1}^{j-r+a} q^{\binom{k-r-a+1}{2}} \qbinom{b-1}{a-1} \qbinom{b+r-a}{k-a} \Theta_{e_{m-j+r}} \Theta_{e_{l+k-j+a}} \nabla E_{j-r+a,b} \\
    \end{align*}
\end{theorem}

Let us fix the following notation:
\begin{equation*}
    \SF(n,k,l) \coloneqq \Theta_{e_k} \Theta_{e_l} \Ht_{(n-k-l)} \rvert_{t=0}
\end{equation*}

We need the following lemma.

\begin{lemma}
    \label{lem:t=0}
    \[ \left( \Theta_{e_{k}} \Theta_{e_{l}} \nabla E_{n-k-l,r} \right)_{t=0} = \delta_{n-k-l, r} \SF(n,k,l). \]
\end{lemma}

\begin{proof}
    By the proof of \cite{IraciRhoadesRomero2023FermionicTheta}*{Lemma~3.6}, we have that, for all $\lambda \vdash n$ and $G \in \Lambda$, \[ \left( \Theta_{e_\lambda} G \right)_{t=0} = \left( \Theta_{e_\lambda} G \rvert_{t=0} \right)_{t=0}. \] In our case, we have that $\nabla E_{n-k-l, r} \rvert_{t=0} = \delta_{n-k-l, r} \Ht_{(n-k-l)}$, so the thesis follows.
\end{proof}

In the rest of this section, we will show that we can rewrite the recursion as follows.

\begin{theorem}
    \label{thm:recursion}
    For any $n,k,l$ with $k+l<n$, $\SF(n,k,l)$ satisfies the recursion
    \begin{align*}
        h_j^\perp & \SF(n,k,l) = \sum_{r=0}^{j} \sum_{a=0}^{j} \sum_{i=0}^{j} \qbinom{n-k-l}{j-r-a+i} q^{\binom{a-i}{2}} \qbinom{n-k-l-(j-r-a+i)}{a-i} \\
                  & \times q^{\binom{r-i}{2}} \qbinom{n-k-l-(j-r-a+i)}{r-i} \qbinom{n-k-l-(j-r-a)-1}{i} \SF(n-j, k-r, l-a)
    \end{align*}
    with initial conditions $\SF(0,k,l) = \delta_{k,0} \delta_{l,0}$, and $\SF(n,k,l) = 0$ if $n<0$.
\end{theorem}

The initial conditions are trivial. In order to remain sane when dealing with all these recursions, we will prove this theorem through a series of intermediate equalities. First, we recall two useful identities.

\begin{proposition}[$q$-Chu-Vandermonde identity]
    \label{prop:q-Chu-Vandermonde}
    \[ \qbinom{j}{a} = \sum_{i=0}^{a} q^{(r-i)(a-i)} \qbinom{r}{i} \qbinom{j-r}{a-i} \]
\end{proposition}

\begin{proposition}[Trinomial identity]
    \label{prop:trinomial}
    \[ \qbinom{x}{y} \qbinom{y}{z} = \qbinom{x}{x-y+z} \qbinom{x-y+z}{z} \]
\end{proposition}

From now on, we use the notation $B \coloneqq n-k-l$ (the number of blocks) to make our formulae more compact.

\begin{lemma}
    \label{lem:intermediate_recursion}
    \begin{align*}
        h_j^\perp \SF(n,k,l) = & \sum_{r=0}^{j} \qbinom{B}{r} \sum_{a=0}^{j} \left( q^{\binom{r-a}{2}} \qbinom{B-(j-r-a)-1}{r-1} \qbinom{j-1}{a} \right. \\
                               & \left. \quad + q^{\binom{r-a+1}{2}} \qbinom{B-(j-r-a)-1}{r} \qbinom{j}{a} \right) \SF(n-j, k-r, l-a).                   \\
    \end{align*}
\end{lemma}

\begin{proof}
    Using \Cref{lem:t=0} and \Cref{thm:full_recursion}, we have
    \begin{align*}
        h_j^\perp \SF(n,k,l)
         & = \sum_{r=0}^{j} \qbinom{B}{r} \sum_{a=0}^{B} \left( q^{\binom{B-r-a}{2}} \qbinom{j-r+a-1}{a} \qbinom{j-1}{B-a-1} \right.   \\
         & \quad + \left. q^{\binom{B-r-a+1}{2}} \qbinom{j-r+a-1}{a-1} \qbinom{j}{B-a} \right) \SF(n-j, k-(j-r), n-k-j-a)              \\
         & = \sum_{r=0}^{j} \qbinom{B}{r} \sum_{a=0}^{B} \left( q^{\binom{a-r}{2}} \qbinom{B+j-r-a-1}{j-r-1} \qbinom{j-1}{a-1} \right. \\
         & \quad + \left. q^{\binom{B-r-a+1}{2}} \qbinom{B+j-r-a-1}{j-r} \qbinom{j}{a} \right) \SF(n-j, k-(j-r), l-(j-a))              \\
         & = \sum_{r=0}^{j} \qbinom{B}{j-r} \sum_{a=0}^{j} \left( q^{\binom{r-a}{2}} \qbinom{B-(j-r-a)-1}{r-1} \qbinom{j-1}{a} \right. \\
         & \quad + \left. q^{\binom{r-a+1}{2}} \qbinom{B-(j-r-a)-1}{r} \qbinom{j}{a} \right) \SF(n-j, k-r, l-a)
    \end{align*}
    as desired. We used $a \mapsto B-a$ in the second equality and $r \mapsto j-r, a \mapsto j-a$ in the last one.
\end{proof}

We are now ready to prove \Cref{thm:recursion}.

\begin{proof}[Proof of \Cref{thm:recursion}]
    First, we apply \Cref{prop:q-Chu-Vandermonde} to \Cref{lem:intermediate_recursion}, getting the following.
    \begin{align*}
        h_j^\perp \SF(n,k,l)
         & = \sum_{r=0}^{j} \qbinom{B}{j-r} \sum_{a=0}^{j} \sum_{i=0}^{a} \left( q^{\binom{r-a}{2}} \qbinom{B-(j-r-a)-1}{r-1} q^{(r-i-1)(a-i)} \qbinom{r-1}{i} \qbinom{j-r}{a-i} \right. \\
         & \left. \quad + q^{\binom{r-a+1}{2}} \qbinom{B-(j-r-a)-1}{r} q^{(r-i)(a-i)} \qbinom{r}{i} \qbinom{j-r}{a-i} \right) \SF(n-j, k-r, l-a).
    \end{align*}
    Now we notice that $\binom{r-a+1}{2} + (r-i)(a-i) = \binom{r-i}{2} + \binom{a-i}{2}$, and we use \Cref{prop:trinomial} with $x=B$, $y=j-r$, and $z=a-i$ (and the corresponding off-by-one term), and get
    \begin{align*}
         & h_j^\perp \SF(n,k,l) = \sum_{r=0}^{j} \sum_{a=0}^{j} \sum_{i=0}^{a} \qbinom{B}{j-r-a+i} q^{\binom{a-i}{2}} \qbinom{B-(j-r-a+i)}{a-i} q^{\binom{r-i}{2}} \\
         & \quad \times \left( \qbinom{B-(j-r-a)-1}{r-1} \qbinom{r-1}{i} + q^{r-i} \qbinom{B-(j-r-a)-1}{r} \qbinom{r}{i} \right) \SF(n-j, k-r, l-a) .
    \end{align*}
    Using again \Cref{prop:trinomial} with $x=B-(j-r-a)-1$, $y=r$, and $z=r-i$ (and the corresponding off-by-one term), we get
    \begin{align*}
        h_j^\perp \SF(n,k,l)
         & = \sum_{r=0}^{j} \sum_{a=0}^{j} \sum_{i=0}^{a} \qbinom{B}{j-r-a+i} q^{\binom{a-i}{2}} \qbinom{B-(j-r-a+i)}{a-i} \\
         & \quad \times q^{\binom{r-i}{2}} \left( \qbinom{B-(j-r-a)-1}{i} \qbinom{B-(j-r-a+i)-1}{r-i-1} \right.            \\
         & \quad \quad \left. + q^{r-i} \qbinom{B-(j-r-a)-1}{i} \qbinom{B-(j-r-a+i)-1}{r-i} \right) \SF(n-j, k-r, l-a)
    \end{align*}
    and finally
    \begin{align*}
        h_j^\perp \SF(n,k,l)
         & = \sum_{r=0}^{j} \sum_{a=0}^{j} \sum_{i=0}^{a} \qbinom{B}{j-r-a+i} q^{\binom{a-i}{2}} \qbinom{B-(j-r-a+i)}{a-i}    \\
         & \quad \times q^{\binom{r-i}{2}} \left( \qbinom{B-(j-r-a+i)-1}{r-i-1} + q^{r-i} \qbinom{B-(j-r-a+i)-1}{r-i} \right) \\
         & \quad \times \qbinom{B-(j-r-a)-1}{i} \SF(n-j, k-r, l-a)                                                            \\
         & = \sum_{r=0}^{j} \sum_{a=0}^{j} \sum_{i=0}^{j} \qbinom{B}{j-r-a+i} q^{\binom{a-i}{2}} \qbinom{B-(j-r-a+i)}{a-i}    \\
         & \quad \times q^{\binom{r-i}{2}} \qbinom{B-(j-r-a+i)}{r-i} \qbinom{B-(j-r-a)-1}{i} \SF(n-j, k-r, l-a)
    \end{align*}
    as desired.
\end{proof}

It is convenient to rewrite the statement as follows.

\begin{corollary}
    \label{cor:recursion}
    Let $\mu \vDash_0 n$ nonzero. Let $m$ be maximal such that $j\coloneqq \mu_m>0$, and define $\mu^- = (\mu_1, \dots, \mu_{m-1},0,0,\dots)$. Then
    \begin{align*}
         & \langle \SF(n,k,l), h_\mu \rangle = \sum_{r=0}^{j} \sum_{a=0}^{j} \sum_{i=0}^{j} \qbinom{n-k-l}{j-r-a+i} q^{\binom{a-i}{2}} \qbinom{n-k-l-(j-r-a+i)}{a-i} \\
         & \quad \times q^{\binom{r-i}{2}} \qbinom{n-k-l-(j-r-a+i)}{r-i} \qbinom{n-k-l-(j-r-a)-1}{i} \langle \SF(n-j,k-r,l-a), h_{\mu^-} \rangle
    \end{align*}
    with initial conditions $\SF(0, k, l) = \delta_{k,0} \delta_{l,0}$.
\end{corollary}

\begin{proof}
    The statement follows immediately from \Cref{thm:recursion}.
\end{proof}

\section{Combinatorial expansion of \texorpdfstring{$\left.\Theta_{e_k}\Theta_{e_l}\nabla e_{n-k-l}\right|_{t=0}$}{the main symmetric function}}
\label{sec:combinatorial_expansion}

Let us define the power series \[ \SW_{x;q}(n,k,l) \coloneqq \sum_{w \in \SW(n,k,l)} q^{\sminv(w)} x_w \] where $x_w = \prod_{i=1}^n x_{w_i}$. Our goal is to show the following theorem, which is our first main result:

\begin{theorem}
    \label{thm:main}
    For any $n,k,l$ with $k+l < n$, we have the identity \[ \SF(n,k,l) =\SW_{x;q}(n,k,l). \]
\end{theorem}

The proof is given in \Cref{sub:general_case}, and is a bit technical. We first give a special case of the result which serves to illustrate some of the combinatorics involved, and is interesting in itself.

\subsection{The standard case}
\label{sub:standard_case}

We take the inner product of the quantity in \Cref{thm:main} with $h_{1^n}$, that is, we want to show $\langle \SF(n,k,l), h_{1^n} \rangle =\SW_{q}(1^n,k,l)$. Recall that the latter is the $q$-enumerator for \emph{segmented permutations} with $k$ ascents and $l$ descents. By specializing \Cref{cor:recursion} to $\mu = 1^n$, so that $j=1$ in that statement, it is readily shown that we have to prove the following recursive formula for $\SW_q(1^n, k, l)$.

\begin{proposition}
    \label{prop:segmented_permutations}
    For any $n,k,l$ with $k+l<n$, the polynomials $\SW_q(1^n,k,l)$ satisfy the recursion
    \begin{align}
        \label{eq:recurrence_segmented_permutations}
        \SW_q(1^n,k,l) = [n-k-l]_q & \left(~\SW_q(1^{n-1},k,l) + \SW_q(1^{n-1},k,l-1)\right. \nonumber \\
                                   & \quad + \left.\SW_q(1^{n-1},k-1,l) +\SW_q(1^{n-1},k-1,l-1)\right).
    \end{align}
    with initial conditions $\SW_q(\varnothing, k, l) = \delta_{k,0} \delta_{l,0}$.
\end{proposition}

Before getting started with the proof, it is worth noticing that in this case the definition of $\sminv$ in \Cref{def:smirnov_inv} simplifies significantly: since there are no equal values, we only have cases~\ref{inv:first_letter} and~\ref{inv:bigger_entry} to consider, which are essentially $2-31$ patterns, cf. \Cref{rem:231}.

\begin{proof}
    Initial conditions are trivial, as the only segmented permutation on $1$ element has no ascent or descent. To get a recursive formula, will insert the value $n$ in a segmented permutation on $n-1$ elements. It can be done in four distinct ways:
    \begin{itemize}
        \item as a new singleton block (keeping the number of ascents and descents the same, and increasing the number of blocks by one), or
        \item at the beginning of a block (creating no ascent and one descent, and keeping the number of blocks the same), or
        \item at the end of a block (creating one ascent and no descents, and keeping the number of blocks the same), or
        \item replacing a block separator (creating both an ascent and a descent, and decreasing the number of blocks by one).
    \end{itemize}

    We want at the end, a permutation $\sigma$ in $\SW(1^n,k,l)$ to prove~\eqref{eq:recurrence_segmented_permutations}. So in these four cases the starting permutation must belong to $\SW(1^{n-1},k,l)$, $\SW(1^{n-1},k,l-1)$, $\SW(1^{n-1},k-1,l)$ and $\SW(1^{n-1},k-1,l-1)$ respectively, and this corresponds to the four terms on the right-hand side of~\Cref{eq:recurrence_segmented_permutations}.

    Note that $\sigma$ has $B\coloneqq n-k-l$ blocks. We claim that each of these four types of insertion can be done in exactly $B$ different ways.

    Indeed, in the first case, the starting segmented permutation must have $B-1$ blocks, so we can insert the new singleton block $n$ in $B$ different positions. In the second (resp. third) case, the starting segmented permutation has $B$ blocks, and we can insert $n$ at the beginning (resp. end) of each of these blocks. In the last case, the starting segmented permutation has $B+1$ blocks, so it has $B$ block separators, and we can replace any of them with $n$.

    In all four cases, we thus have $B$ possibilities of insertion, and $n$ forms a sminversion with all and only the initial positions to its right, thus contributing the factor $1+q+\cdots+q^{B-1}=[B]_q$. Moreover, since $n$ is the biggest entry so far, all the existing sminversions remain sminversions (suitably shifted by the position of $n$) as is easily checked. This proves the recursive formula~\eqref{eq:recurrence_segmented_permutations}.
\end{proof}

\Cref{prop:segmented_permutations} is of interest in itself because it provides a combinatorial formula for the polynomial corresponding to the conjectured \emph{Hilbert series} of the $(1,2)$-coinvariant module, namely \[ \sum_{k+l<n} \langle \Theta_{e_k} \Theta_{e_l} \nabla e_{(n-k-l)}, h_{1^n} \rangle \rvert_{t=0} u^k v^l. \]
We now tackle the general case, which gives the full conjectured Frobenius characteristic of the same module.

\subsection{The general case}
\label{sub:general_case}

In order to better understand the following proof, the reader is invited to check \Cref{ex:recursion} and \Cref{ex:recursion_big} while reading it. As in the standard case above,  we have four kinds of insertions: the difference is that now the maximal letter can be inserted multiple times, and we have to be a little careful about the order in which we perform the insertions. We will detail the approach here because it is going to be convenient when we come back to it in \Cref{sec:link-delta}.

The idea is the following: let $\mu \vDash_0 n$, $m$ maximal such that $j \coloneqq \mu_m > 0$, and define $\mu^- = (\mu_1, \dots, \mu_{m-1},0,\dots)$. We can build a word $\overline{w} \in \SW(\mu)$ from a word $w \in \SW(\mu^-)$ by inserting $j$ occurrences of $m$ in $w$, which contains only letters smaller than $m$.

There are four kinds of insertions, depending on whether $m$ is a peak (replacing a block separator), double fall (at the beginning of a block), double rise (at the end of a block), or valley (in a new, singleton block) in $\overline{w}$; for each kind of insertion, we have to show that sminversions of $w$ remain sminversions in $\overline{w}$, and we have to keep track of the amount of new sminversions the occurrences of $m$ form with the remaining letters of $\overline{w}$. The way we do so is by ordering the possible insertion positions from right to left, corresponding to increasing contributions, and then computing the relevant $q$-enumerators.

Let us start with a preliminary lemma.

\begin{lemma}
    \label{lem:preserve_inversions}
    Let $w \in \SW(n)$, and let $m \geq \max w$. Let $\overline{w} \in \SW(n+1)$ be a word obtained from $w$ by inserting an occurrence of $m$ as a double fall, double rise, or singleton. Then, if $(i,j)$ forms a sminversion in $w$, the corresponding indices in $\overline{w}$ also form a sminversion. If instead $m$ is inserted as a peak, the same holds as long as $m > \max w$.
\end{lemma}

\begin{proof}
    Suppose that $(i,j)$ is a sminversion in $w$; it is clear that, if $m$ is inserted anywhere except in position $j$, then the corresponding pair remains a sminversion. Suppose now that $m$ is inserted in position $j$: we want to show that $(i,j+1)$ is a sminversion of $\overline{w}$. If $w_i < m$, now we have $w_i < w_j$ and $w_i > w_{j+1}$, and so $(i,j+1)$ forms a sminversion. If $w_i = m$, then the newly inserted $m$ cannot be a peak (because $(i,j)$ formed a sminversion in $w$ and $m = \max w$), and so either $m$ is a valley or a double rise, in which case $j+1$ is initial (as in \Cref{def:smirnov_inv}, case~\ref{inv:first_letter}), or $m$ is inserted as a double fall, in which case it is must be initial, and so $(i,j+1)$ is a sminversion of the kind \Cref{def:smirnov_inv}, case~\ref{inv:equal_then_first_letter}.
\end{proof}

It is apparent from \Cref{lem:preserve_inversions} that we should insert peaks first, as doing so when there are already occurrences of $m$ in $w$ can affect the number of sminversions.

\begin{lemma}
    \label{lem:peak_insertion}
    Let $w \in \SW(n,k,l)$, and let $m > \max w$. The $q$-enumerator with respect to $\sminv$ of words obtained from $w$ by inserting $s$ occurrences of $m$ in $w$ as a peak (i.e.\ replacing a block separator) is \[ \qbinom{n-k-l-1}{s} q^{\sminv(w)}. \]
\end{lemma}

\begin{proof}
    Let $\overline{w}$ be a word obtained as in the statement. Since $m > \max w$, conditions \ref{inv:bigger_entry}, \ref{inv:equal_then_first_letter} and \ref{inv:equal_then_bigger_entry} in \Cref{def:smirnov_inv} never happen, and so each occurrence of $m$ forms a sminversion with each and only the initial letters of the blocks of $\overline{w}$ to its right. Since $w$ has $n-k-l-1$ block separators, or equivalently $\overline{w}$ has $n-k-l-s$ blocks, these sminversions are $q$-counted by $\qbinom{n-k-l-1}{s}$, which concludes the proof.
\end{proof}

\begin{lemma}
    \label{lem:fall_insertion}
    Let $w \in \SW(n,k,l)$, and let $m \geq \max w$, with no initial $m$ in $w$. The $q$-enumerator with respect to $\sminv$ of words obtained from $w$ by inserting $s$ occurrences of $m$ in $w$ as a double fall (i.e.\ as an initial element) is \[ q^{\binom{s}{2}} \qbinom{n-k-l}{s} q^{\sminv(w)}. \]
\end{lemma}

\begin{proof}
    As in \Cref{lem:peak_insertion}, each occurrence of $m$ forms a sminversion with each and only the initial letters of the blocks of $w$ (or $\overline{w}$, as they have the same blocks) to its right. Unlike \Cref{lem:peak_insertion}, however, we can insert at most one occurrence of $m$ in each block. Since $\overline{w}$ has $n-k-l$ blocks, the new sminversions are $q$-counted by $q^{\binom{s}{2}} \qbinom{n-k-l}{s}$, which concludes the proof.
\end{proof}

\begin{lemma}
    \label{lem:rise_insertion}
    Let $w \in \SW(n,k,l)$, and let $m \geq \max w$ with no final $m$ in $w$. The $q$-enumerator with respect to $\sminv$ of words obtained from $w$ by inserting $s$ occurrences of $m$ in $w$ as a double rise (i.e.\ as the final element of a block) is \[ q^{\binom{s}{2}} \qbinom{n-k-l}{s} q^{\sminv(w)}. \]
\end{lemma}

\begin{proof}
    As in \Cref{lem:fall_insertion}, each occurrence of $m$ forms a sminversion with each and only the initial letters of the blocks of $\overline{w}$ to its right, and we can insert at most one occurrence of $m$ in each block. Since $w$ has $n-k-l$ blocks, the new sminversions are $q$-counted by $q^{\binom{s}{2}} \qbinom{n-k-l}{s}$, which concludes the proof.
\end{proof}

\begin{lemma}
    \label{lem:valley_insertion}
    Let $w \in \SW(n,k,l)$, and let $m \geq \max w$ with no singleton $m$ in $w$. The $q$-enumerator with respect to $\sminv$ of words obtained from $w$ by inserting $s$ occurrences of $m$ in $w$ as a valley (i.e.\ as a singleton block) is \[ \qbinom{n-k-l+s}{s} q^{\sminv(w)}. \]
\end{lemma}

\begin{proof}
    Each occurrence of $m$ creates exactly one sminversion with each block of $w$ to its right. Indeed, since $m$ is maximal, conditions \ref{inv:bigger_entry} and \ref{inv:equal_then_bigger_entry} in \Cref{def:smirnov_inv} never occur; if a block to the right of our occurrence of $m$ begins with an entry that is strictly less than $m$, then condition \ref{inv:first_letter} occurs and we have a sminversion; if it begins with an $m$, then the following letter (which exists because it is not a singleton) is strictly smaller than $m$, so condition \ref{inv:equal_then_first_letter} occurs. Since $w$ has $n-k-l$ blocks, we have $n-k-l+1$ positions where we can insert an $m$, each creating a number of new sminversions equal to the number of blocks of $w$ to its right, for a total contribution of $\qbinom{n-k-l+s}{s}$, as desired.
\end{proof}

We are now ready to put the pieces together.

\begin{proof}[Proof of \Cref{thm:main}]
    We prove the equivalent statement \[ \langle \SF(n,k,l), h_\mu \rangle = \SW_q(\mu, k, l). \]

    Our goal is to show that the right hand side satisfies the same recurrence as the left hand side, that is the one given in \Cref{cor:recursion}. In other words, we have to show that, for any nonzero $\mu \vDash_0 n$, with $m$ maximal such that $j\coloneqq \mu_m>0$, and $\mu^- = (\mu_1, \dots, \mu_{m-1},0,\dots)$, we have
    \begin{align}
        \label{eq:recurrenceSW}
         & \SW_q (\mu, k, l) = \sum_{r=0}^{j} \sum_{a=0}^{j} \sum_{i=0}^{j} \qbinom{n-k-l}{j-r-a+i} q^{\binom{a-i}{2}} \qbinom{n-k-l-(j-r-a+i)}{a-i} \\
         & \quad \times q^{\binom{r-i}{2}} \qbinom{n-k-l-(j-r-a+i)}{r-i} \qbinom{n-k-l-(j-r-a)-1}{i} \SW_q(\mu^-, k-r, l-a),
    \end{align}
    with initial conditions $\SW_{q}(\varnothing,k,l) = \delta_{k,0} \delta_{l,0}$.

    It is clear that initial conditions are satisfied, as the only word of length $0$ has no ascents or descents. If $w \in \SW_q(\mu, k, l)$, then $m$ is the greatest letter appearing in $w$, and it appears exactly $j$ times. We interpret the variables appearing in the recursion as follows:
    \begin{itemize}
        \item $i$ is the number of occurrences of $m$ that are neither at the beginning nor at the end of a block;
        \item $a-i$ is the number of occurrences of $m$ that are at the beginning of a block of size at least $2$;
        \item $r-i$ is the number of occurrences of $m$ that are at the end of a block of size at least $2$;
        \item $j-r-a+i$ is the number of occurrences of $m$ that are singletons.
    \end{itemize}

    We proceed backwards, starting from a word $w \in \SW_{q}(\mu^-, k-r, l-a)$ for all the possible values of $r$ and $a$, and inserting $j$ occurrences of $m$.

    First, for each possible value of $i$, we insert $i$ peaks. By \Cref{lem:peak_insertion}, since $w$ has $n-k-l-(j-r-a)$ blocks, this yields a contribution of $\qbinom{n-j-(k-r)-(l-a)-1}{i}$, and we are left with a word with $n-k-l-(j-r-a+i)$ blocks, $k-(r-i)$ ascents, and $l-(a-i)$ descents.

    Next, we insert $a-i$ double falls, and $r-i$ double rises. By \Cref{lem:fall_insertion} and \Cref{lem:rise_insertion}, the former yields a contribution of $q^{\binom{r-i}{2}} \qbinom{n-k-l-(j-r-a+i)}{r-i}$, and the latter yields a contribution of $q^{\binom{a-i}{2}} \qbinom{n-k-l-(j-r-a+i)}{a-i}$. We are left with a word with the same number of blocks, $k$ ascents, and $l$ descents.

    Finally, we insert $j-r-a+i$ singletons. By \Cref{lem:valley_insertion}, this yields a contribution of $\qbinom{n-k-l}{j-r-a+i}$, and we are left with a word in $\SW(n,k,l)$, as expected. By construction, each of these words appears exactly once in this process, so the thesis follows.
\end{proof}

We now show two examples for this process: in \Cref{ex:recursion} we start from a small word and go through all the possible insertions; in \Cref{ex:recursion_big} we start from a bigger word and choose one possible insertion.

\begin{example}
    \label{ex:recursion}
    We take $\mu = (5,3,6,0,\ldots)$ so that $j = 6$ and $\mu^-=(5,3,0,\ldots)$ and $k = 6$, $l = 5$. Furthermore, we fix $r = 4$, $a = 3$ and $i=2$, and we will study the summand of the recursion for these values.

    We start by picking an example of an element in $\SW(\mu^-, k-r, l-a)= \SW((5,3), 2, 2)$ \[w_0 \coloneqq 1|21|12|121.\]
    The $\sminv$ of this word is equal to $3$. From this word, we will build some elements of $\SW(\mu,k,l) = \SW((5,3,6),6,5)$ that reduce to $w_0$ by deleting its maximal letters that are the first or last letters of its blocks, deleting singleton blocks containing a maximal letter, and replacing all other maximal letters with a block separator. We will keep track of the number of sminversions during the process.
    \begin{itemize}
        \item Since we took $i = 2$, we must insert $2$ maximal letters that are neither the first, nor the last letter of its block. The number of block separators of an element in $\SW(\mu^-, k-r, l-a)$ is equal to $n-j-(k-r)-(l-a)-1=3$. Pick $i=2$ among these block separators, and replace them with a maximal letter:
              \[
                  \begin{forest}
                      [1|21|12|121
                      [1321312|121, edge label = {node[pos = .2,left=5pt, font = \scriptsize]{+2}}],
                      [1321|123121, edge label = {node[pos = .5,left, font = \scriptsize]{+1}}],
                      [1|213123121, edge label = {node[pos = .2,right=5pt, font = \scriptsize]{+0}}]
                      ]
                  \end{forest}
              \]
              and so the extra sminversions are accounted for by the factor \[\qbinom{n-k-(k-r)-(l-a)-1}{i} = \qbinom{3}{2} = q^2 + q + 1.\]
              Let us continue our example with the word $1321312|121$ of $\sminv$ equal to $3+2 = 5$.
        \item Since $r-i = 2$, we must insert $2$ maximal letters at the beginning of a block. Since an element of $\SW(\mu^-, k-r, l-a)$ has $n-j-(k-r)-(l-a)=4$ blocks, and we removed $i=2$ blocks at the previous step, our words now have $n-k-l - (j-r-a+i) = 2$ blocks, of which we must choose $a-i=2$.
              \[
                  \begin{forest}
                      [1321312|121
                      [31321312|3121, edge label = {node[pos = .5,left, font = \scriptsize]{+1}}]
                      ]
                  \end{forest}
              \]
              Doing this in the only way possible we get the expected number of sminversions.
              \[q^{\binom{a-i}{2}}\qbinom{n-k-l - (j-r-a+i)}{a-i} = q^{\binom{2}{2}}\qbinom{2}{2} = q.\]
              We continue our example with the word $31321312|3121$ of $\sminv$ equal to $5+1 = 6$
        \item Now we insert $r-i = 1$ maximal letter at the end of blocks.
              \[
                  \begin{forest}
                      [31321312|3121
                      [313213123|3121, edge label = {node[pos = .5,left, font = \scriptsize]{+1}}]
                      [31321312|31213, edge label = {node[pos = .5,right, font = \scriptsize]{+0}}]
                      ]
                  \end{forest}
              \]
              A similar argument as before gets us the right number of sminversions
              \[q^{\binom{r-i}{2}}\qbinom{n-k-l - (j-r-a+i)}{r-i} = q^{\binom{1}{2}}\qbinom{2}{1} = 1+q.\]
              Let us continue our example with the word $313213123|3121$ of $\sminv$ equal to $6+1=7$.
        \item Finally, we insert $j-r-a+i = 1$ maximal letter in singleton blocks.
              \[
                  \begin{forest}
                      [313213123|3121
                      [3|313213123|3121, edge label = {node[pos = .2,left = 8pt, font = \scriptsize]{+2}}]
                      [313213123|3|3121, edge label = {node[pos = .5,left, font = \scriptsize]{+1}}]
                      [313213123|3121|3, edge label = {node[pos = .2,right = 8pt, font = \scriptsize]{+0}}]
                      ]
                  \end{forest}
              \]
              We had to interlace the $n-k-l-(j-r-a+i)$ existing blocks with the $j-r-a+i$ new singleton blocks, giving a $\sminv$ contribution of \[\qbinom{n-k-l}{j-r-a+i} = \qbinom{3}{1} = 1 + q + q^2.\] These last three words all belong to $\SW(\mu,k,l) = \SW((5,3,6),6,5)$ and have $\sminv$ equal to $7+2$, $7+1$, and $7+0$, respectively.
    \end{itemize}
\end{example}

\begin{example}
    \label{ex:recursion_big}
    Consider the word
    \[w_0=23\underset{\circ}{|}121\underset{\circ}{|}32\underset{\bullet}{|}1\underset{\circ}{|}1231\underset{\bullet}{|}12\underset{\circ}{|}23,\] where the disks indicate possible peak insertion positions, and the filled ones corresponding to the chosen ones. After these insertions, we obtain
    \[w_1={}\underset{\bullet}{}23{|\underset{\circ}{}}121{|\underset{\bullet}{}}32{\bf 4}1{|\underset{\circ}{}}1231{\bf 4}12{|\underset{\bullet}{}}23.\]
    Insertion positions of potential and chosen initial elements are illustrated, leading to
    \[w_2={\bf 4}23\underset{\circ}{}{|}121\underset{\bullet}{}{|}{\bf 4}32{4}1\underset{\bullet}{}{|}1231{4}12\underset{\circ}{}{|}{\bf 4}23\underset{\circ}{}.\]
    Now insertion positions of potential and actual final elements are illustrated, and we have
    \[w_3=\underset{\circ_1}{}{4}23\underset{\circ_0}{|}121{\bf 4}\underset{\circ_0}{|}{4}32{4}1{\bf 4}\underset{\circ_3}{|}1231{4}12\underset{\circ_1}{|}{4}23\underset{\circ_0}{}.\]
    There possible insertion positions of singleton blocks are illustrated, the index referring to the number of such blocks to be inserted there, leading us to the final word
    \[w_4={\bf 4}|{4}23{|}121{4}{|}{4}32{4}1{4}{|}{\bf 4}|{\bf 4}|{\bf 4}{|}1231{4}12{|}{\bf 4}{|}{4}23.\]
    In all, we successively inserted $2,2,3$ and $5$ letters corresponding to the values $i=2,a=5,r=4$ and $j=12$ in the proof. This led from a word $w_0$ in $\SW(16,6,3)$ to a word $w_4$ in $\SW(28,10,8)$.

    As for sminversions, we added $(1+3)$ sminversions from $w_0$ to $w_1$, $3+(0+2+4)$ sminversions from $w_1$ to $w_2$, $1+(2+3)$ sminversions from $w_2$ to $w_3$, and finally $1+3+4+5+9$ sminversions from $w_3$ to $w_4$. The corresponding $q$-enumerators as in Lemmas~\ref{lem:peak_insertion} to~\ref{lem:valley_insertion} are successively
    \[\qbinom{6}{2}q^{\sminv{(w_0)}},\quad q^3\qbinom{5}{3}q^{\sminv{(w_1)}},\quad q\qbinom{5}{2}q^{\sminv{(w_2)}},\quad \qbinom{11}{5}q^{\sminv{(w_3)}}.\]
\end{example}

\section{Unified Delta theorem at \texorpdfstring{$t=0$}{t=0}}
\label{sec:link-delta}

Theorem~\ref{thm:main} gives us a combinatorial interpretation of $\SF(n,k,l)$ in terms of segmented Smirnov words. Recall that \[ \SF(n,k,l) \coloneqq \Theta_{e_k} \Theta_{e_l} \Ht_{(n-k-l)} \rvert_{t=0} = \Theta_{e_k} \Theta_{e_l} \nabla e_{n-k-l} \rvert_{t=0},\] the latter equality being \cite{IraciRhoadesRomero2023FermionicTheta}*{Lemma~3.6} (it also follows from \Cref{lem:t=0} by taking a sum over $r$). We recognize this as the symmetric function of the unified Delta Conjecture \eqref{eq:unified-delta}, evaluated at $t=0$. In this section, we interpret a variant of the combinatorics of \Cref{thm:main} in the setting of labelled Dyck paths, which was the original motivation. As we will see, our result reinforces the idea that there should be a $q$-statistic on labelled Dyck paths with both decorated rises and decorated contractible valleys that extends both the rise and the valley version of the Delta conjecture.

We will define a $\dinv$ statistic on decorated Dyck paths of area zero and show the following.

\begin{theorem}[Unified Delta Theorem at $t=0$]
    \label{thm:unified_Delta_t_equals_0}
    \[ \left. \Theta_{e_k} \Theta_{e_l} \nabla e_{n-k-l} \right\rvert_{t=0} = \sum_{D \in\LD_0(n)^{\ast k, \bullet l}} q^{\dinv(D)} x^D. \]
\end{theorem}

\subsection{The Delta conjecture}
\label{sub:delta_conjecture}

It is now time to give some more precise definitions of the combinatorics of the Delta conjecture.

\begin{definition}
    A \emph{Dyck path} of size $n$ is a lattice path starting at $(0,0)$, ending at $(n,n)$, using only unit North and East steps, and staying weakly above the line $x=y$. A \emph{labelled Dyck path} is a Dyck path together with a positive integer label on each of its vertical steps such that labels on consecutive vertical steps must be strictly increasing (from bottom to top). We will draw the labels of vertical steps in the square to its right.

    A \emph{rise} of a labelled Dyck path is a vertical step that is preceded by another vertical step.

    A \emph{valley} of a labelled Dyck path is a vertical step $v$ preceded by a horizontal step. A valley $v$ is \emph{contractible} if it is either preceded by two horizontal steps, or by a horizontal step that is itself preceded by a vertical step whose label is strictly smaller than $v$'s label.

    A \emph{decorated labelled Dyck path} is a labelled Dyck path, together with a choice of rises and contractible valleys, which are \emph{decorated}.
    We set
    \begin{align*}
         & \mathsf{DRise}(D) = \{i \in [n] \mid \text{the $i\th$ vertical step of $D$ is a decorated rise}\}     \\
         & \mathsf{DValley}(D) = \{i \in [n] \mid \text{the $i\th$ vertical step of $D$ is a decorated valley}\}
    \end{align*}
    We decorate rises with a $\ast$ and valleys with a $\bullet$, and these decorations will be displayed in the square to the left of the vertical step. The set of decorated labelled Dyck paths of size $n$ with $k$ decorated rises and $l$ decorated valleys, is denoted by $\LD(n)^{\ast k, \bullet l}$.
\end{definition}

\begin{figure}
    \centering
    \begin{minipage}{.48\textwidth}
        \begin{tikzpicture}[scale=.72]
            \draw[step=1.0, gray!60, thin] (0,0) grid (8,8);
            \draw[gray!60, thin] (0,0) -- (8,8);
            \draw[blue!60, line width=3pt] (0,0) -- (0,1) -- (0,2) -- (1,2) -- (1,3) -- (2,3) -- (2,4) -- (2,5) -- (2,6) -- (3,6) -- (4,6) -- (4,7) -- (5,7) -- (6,7) -- (7,7) -- (7,8) -- (8,8);
            \node at (0.5,0.5) {$2$};
            \draw (0.5,0.5) circle (.4cm);
            \node at (0.5,1.5) {$3$};
            \draw (0.5,1.5) circle (.4cm);
            \node at (1.5,2.5) {$4$};
            \draw (1.5,2.5) circle (.4cm);
            \node at (2.5,3.5) {$1$};
            \draw (2.5,3.5) circle (.4cm);
            \node at (2.5,4.5) {$2$};
            \draw (2.5,4.5) circle (.4cm);
            \node at (2.5,5.5) {$4$};
            \draw (2.5,5.5) circle (.4cm);
            \node at (4.5,6.5) {$3$};
            \draw (4.5,6.5) circle (.4cm);
            \node at (7.5,7.5) {$2$};
            \draw (7.5,7.5) circle (.4cm);
            \node at (1-1-0.5,1+0.5) {$\ast$};
            \node at (5-3-0.5,5+0.5) {$\ast$};
            \node at (2-1-0.5,2+0.5) {$\bullet$};
            \node at (6-2-0.5,6+0.5) {$\bullet$};
        \end{tikzpicture}
    \end{minipage}%
    \begin{minipage}{.48\textwidth}
        \begin{tikzpicture}[scale = .72]
            \draw[draw=none, use as bounding box] (-1, -1) rectangle (9,9);
            \draw[step=1.0, gray!60, thin] (0,0) grid (8,8);

            \draw[gray!60, thin] (0,0) -- (8,8);

            \draw[blue!60, line width=3pt] (0,0) -- (0,1) -- (0,2) -- (0,3) -- (1,3) -- (2,3) -- (3,3) -- (3,4) -- (3,5) -- (4,5) -- (5,5) -- (5,6) -- (6,6) -- (6,7) -- (6,8) -- (7,8) -- (8,8);

            \node at (0.5,0.5) {$1$};
            \draw (0.5,0.5) circle (.4cm);
            \node at (0.5,1.5) {$2$};
            \draw (0.5,1.5) circle (.4cm);
            \node at (0.5,2.5) {$4$};
            \draw (0.5,2.5) circle (.4cm);
            \node at (3.5,3.5) {$2$};
            \draw (3.5,3.5) circle (.4cm);
            \node at (3.5,4.5) {$4$};
            \draw (3.5,4.5) circle (.4cm);
            \node at (5.5,5.5) {$1$};
            \draw (5.5,5.5) circle (.4cm);
            \node at (6.5,6.5) {$3$};
            \draw (6.5,6.5) circle (.4cm);
            \node at (6.5,7.5) {$4$};
            \draw (6.5,7.5) circle (.4cm);
            \node at (1-1-0.5,1+0.5) {$\ast$};
            \node at (2-2-0.5,2+0.5) {$\ast$};
            \node at (4-1-0.5,4+0.5) {$\ast$};
            \node at (7-1-0.5,7+0.5) {$\ast$};
            \node at (3-0-0.5,3+0.5) {$\bullet$};
            \node at (6-0-0.5,6+0.5) {$\bullet$};
        \end{tikzpicture}
    \end{minipage}
    \caption{An element of $\LD(8)^{\ast 2, \bullet 2}$ (left) and an element of $\LD_0(8)^{\ast 4,\bullet 2}$ (right)}
    \label{fig:path-example}
\end{figure}
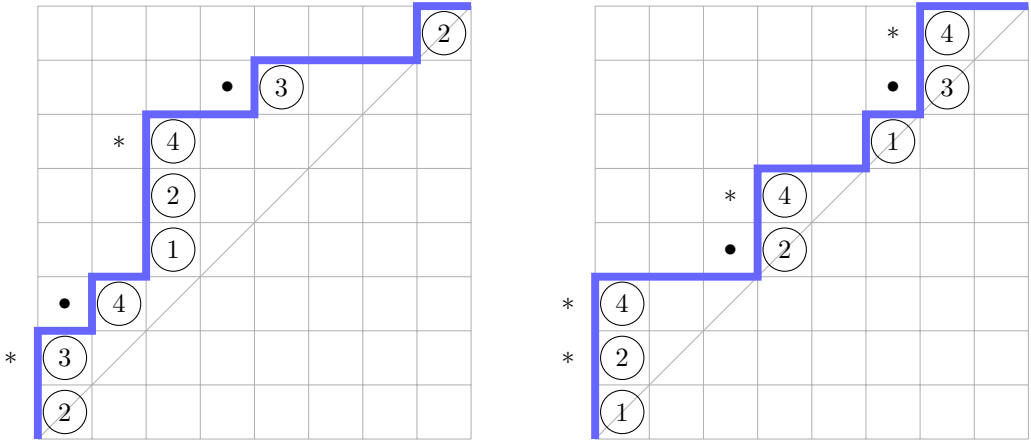

See \Cref{fig:path-example} for an example of an element in $\LD(8)^{\ast 2, \bullet 2}$.

\begin{definition}[Area]
    \label{def:area}
    Given a decorated labelled Dyck path $D$ of size $n$, its \emph{area word} is the word of non-negative integers whose $i\th$ letter equals the number of whole squares between the $i\th$ vertical step of the path and the line $x=y$. If $a$ is the area word of $D$, the \emph{area} of $D$ is \[\area(D) \coloneqq \sum_{i\in [n]\setminus \mathsf{DRise}(D)} a_i.\]
\end{definition}

\begin{definition}[Diagonal inversions]
    \label{def:dinv}
    Take $D$ a decorated labelled Dyck path $D$ of size $n$ and area word $a$. Set $l$ to be the word such that $l_i$ is the label of the $i\th$ vertical step of $D$. For $1\leq i < j \leq n$, we say that
    \begin{itemize}
        \item $(i,j)$ is a \emph{primary diagonal inversion} if $a_i = a_j$, $l_i < l_j$ and $i\not\in \mathsf{DValley}(D)$;
        \item $(i,j)$ is a \emph{secondary diagonal inversion} if $a_i = a_j + 1$, $l_i > l_j$ and $i\not\in \mathsf{DValley}(D)$.
    \end{itemize} Then we define $\displaystyle{\dinv(D)\coloneqq \#\{1\leq i < j \leq n \mid \text{$(i,j)$ is a diagonal inversion}\} - \# \mathsf{DValley}(D)}$.
\end{definition}

\begin{example}
    Take $D$ to be the path on the left in \Cref{fig:path-example}. We have $\mathsf{DRise}(D)=\{2,6\}$, $\mathsf{DValley} = \{3,7\}$. The area word of $D$ is $0 1 1 1 2 3 2 0$, and so its area equals 6. There are $2$ primary diagonal inversions: $(2,3)$ and $(5,7)$; and $2$ secondary diagonal inversions $(2,8)$ and $(6,7)$. Since $\#\mathsf{DValley}(D) = 2$, then $\dinv$ is equal to $4-2=2$.
\end{example}

For $\mu \vDash_0 n$ let $\LD(\mu)^{\ast k, \bullet l}$ be the subset of Dyck paths in $\LD(n)^{\ast k, \bullet l}$ with content $\mu$, that is, their multiset of labels is $\{i^{\mu_i} \mid i > 0\}$. For $D \in \LD(\mu)^{\ast k, \bullet l}$, let $x^D \coloneqq x^\mu$. Then the Delta conjecture states that \[\Theta_{e_k} \nabla e_{n-k} = \sum_{D \in \LD(n)^{\ast k}} q^{\dinv(D)} t^{\area(D)} x^D = \sum_{D \in \LD(n)^{\bullet k}} q^{\dinv(D)} t^{\area(D)} x^D, \] where the first equality (the \emph{rise version}) is now a theorem \cite{DAdderioMellit2022CompositionalDelta}.

\subsection{From segmented Smirnov words to decorated Dyck paths}
\label{sub:bijection_words_to_paths}

We have shown in \Cref{prop:sminv_to_inv} that, when $k=0$ or $l=0$, segmented Smirnov words reduce to ordered set partitions. Combining results from various papers in the literature \cites{RemmelWilson2015ExtensionMacmahonsEquidistribution, Wilson2016Equidistribution, Rhoades2018OrderedSetPartition, HaglundRemmelWilson2018DeltaConjecture}, it has been shown that ordered set partitions have several equidistributed statistics, namely $\inv$, $\dinv$, $\maj$ and $\minimaj$, which in turn have been shown to bijectively match all the specializations at $q=0$ or $t=0$ of the various versions of the Delta conjecture.

For reasons that will be clear in a moment, we are mainly interested in the $\inv$ and $\dinv$ statistics. We recalled the definition of the former in \Cref{def:op_inv}; the latter is as follows.

\begin{definition}
    \label{def:op_dinv}
    For an ordered set partition $\pi = (\pi_1, \dots, \pi_r) \in \OP(\mu, r)$, a diagonal inversion is a triple $h, i, j$ such that either $i < j$ and the $h\th$ smallest element in $\pi_i$ is strictly greater than the $h\th$ smallest element in $\pi_j$, or $i > j$ and the $(h+1)\th$ smallest element in $\pi_i$ is strictly greater than the $h\th$ smallest element in $\pi_j$. We denote by $\dinv(\pi)$ the number of diagonal inversions of $\pi$.
\end{definition}

Now let $\LD_0(n)^{\ast k, \bullet l}$ be the subset of area $0$ Dyck paths in $\LD(n)^{\ast k, \bullet l}$, that is the subset of paths where every rise is decorated and every valley lies on the main diagonal; see for example the path on the right in \Cref{fig:path-example}. Also let $\LD_0(\mu)^{\ast k, \bullet l} \coloneqq \LD(\mu)^{\ast k, \bullet l} \cap \LD_0(n)^{\ast k, \bullet l}$.

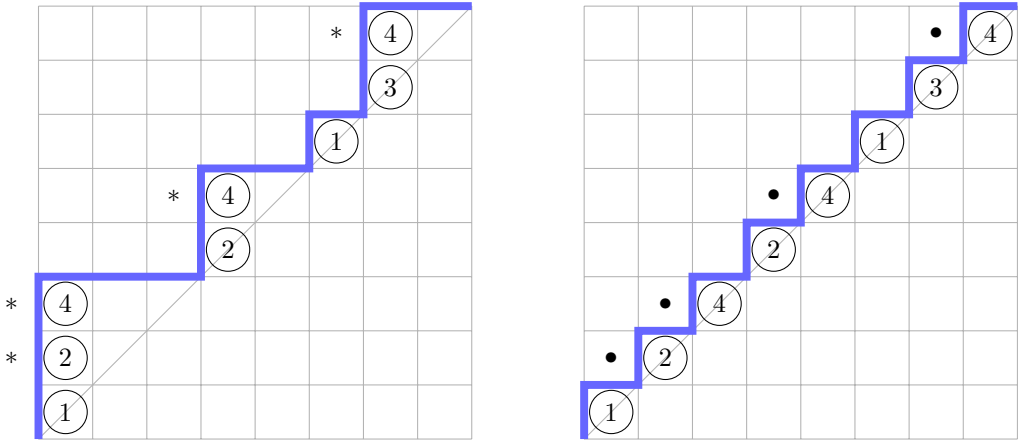
\begin{figure}
    \centering
    \begin{minipage}{.48\textwidth}
        \begin{tikzpicture}[scale = .72]
            \draw[draw=none, use as bounding box] (-1, -1) rectangle (9,9);
            \draw[step=1.0, gray!60, thin] (0,0) grid (8,8);

            \draw[gray!60, thin] (0,0) -- (8,8);

            \draw[blue!60, line width=3pt] (0,0) -- (0,1) -- (0,2) -- (0,3) -- (1,3) -- (2,3) -- (3,3) -- (3,4) -- (3,5) -- (4,5) -- (5,5) -- (5,6) -- (6,6) -- (6,7) -- (6,8) -- (7,8) -- (8,8);

            \node at (0.5,0.5) {$1$};
            \draw (0.5,0.5) circle (.4cm);
            \node at (0.5,1.5) {$2$};
            \draw (0.5,1.5) circle (.4cm);
            \node at (0.5,2.5) {$4$};
            \draw (0.5,2.5) circle (.4cm);
            \node at (3.5,3.5) {$2$};
            \draw (3.5,3.5) circle (.4cm);
            \node at (3.5,4.5) {$4$};
            \draw (3.5,4.5) circle (.4cm);
            \node at (5.5,5.5) {$1$};
            \draw (5.5,5.5) circle (.4cm);
            \node at (6.5,6.5) {$3$};
            \draw (6.5,6.5) circle (.4cm);
            \node at (6.5,7.5) {$4$};
            \draw (6.5,7.5) circle (.4cm);
            \node at (1-1-0.5,1+0.5) {$\ast$};
            \node at (2-2-0.5,2+0.5) {$\ast$};
            \node at (4-1-0.5,4+0.5) {$\ast$};
            \node at (7-1-0.5,7+0.5) {$\ast$};
        \end{tikzpicture}
    \end{minipage}%
    \begin{minipage}{.48\textwidth}
        \begin{tikzpicture}[scale = .72]
            \draw[draw=none, use as bounding box] (-1, -1) rectangle (9,9);
            \draw[step=1.0, gray!60, thin] (0,0) grid (8,8);

            \draw[gray!60, thin] (0,0) -- (8,8);

            \draw[blue!60, line width=3pt] (0,0) -- (0,1) -- (1,1) -- (1,2) -- (2,2) -- (2,3) -- (3,3) -- (3,4) -- (4,4) -- (4,5) -- (5,5) -- (5,6) -- (6,6) -- (6,7) -- (7,7) -- (7,8) -- (8,8);

            \node at (0.5,0.5) {$1$};
            \draw (0.5,0.5) circle (.4cm);
            \node at (1.5,1.5) {$2$};
            \draw (1.5,1.5) circle (.4cm);
            \node at (2.5,2.5) {$4$};
            \draw (2.5,2.5) circle (.4cm);
            \node at (3.5,3.5) {$2$};
            \draw (3.5,3.5) circle (.4cm);
            \node at (4.5,4.5) {$4$};
            \draw (4.5,4.5) circle (.4cm);
            \node at (5.5,5.5) {$1$};
            \draw (5.5,5.5) circle (.4cm);
            \node at (6.5,6.5) {$3$};
            \draw (6.5,6.5) circle (.4cm);
            \node at (7.5,7.5) {$4$};
            \draw (7.5,7.5) circle (.4cm);
            \node at (0.5,1+0.5) {$\bullet$};
            \node at (1.5,2+0.5) {$\bullet$};
            \node at (3.5,4+0.5) {$\bullet$};
            \node at (6.5,7+0.5) {$\bullet$};
        \end{tikzpicture}
    \end{minipage}
    \caption{The images of the ordered set partition $34|1|24|124$ (or $43|1|42|421$) via the two bijections.}
    \label{fig:op-bijection}
\end{figure}

We care in particular about the bijections
\begin{equation}
    \label{eq:dinv_bijection} \OP(\mu, n-k) \leftrightarrow \LD_0(\mu)^{\ast k}
\end{equation}
\begin{equation}
    \label{eq:inv_bijection} \OP(\mu, n-l) \leftrightarrow \LD_0(\mu)^{\bullet l}
\end{equation}
respectively mapping the $\dinv$ and $\inv$ statistics on ordered set partitions to the $\dinv$ on Dyck paths. See \Cref{fig:op-bijection} for an example or \cite{HaglundRemmelWilson2018DeltaConjecture}*{Proposition~4.1} for more details.

\begin{remark}
    While not relevant to our goal, we would like to point out that \eqref{eq:dinv_bijection} also sends $\minimaj$ to $\pmaj$, thus proving the $q=0$ specialization of the $\pmaj$ version of the Delta conjecture \cite{DAdderioIraciVandenWyngaerd2022TheBible}*{Conjecture~2.5}; while the proof of this fact is trivial, the authors are not aware of it being mentioned anywhere in the literature. The case $t=0$ of the same conjecture remains open.
\end{remark}

Using again an insertion method, we will show a bijection between segmented Smirnov words of size $n$ with $k$ ascents and $l$ descents, and labelled Dyck paths of size $n$ with $k$ decorated rises and $l$ decorated contractible valleys, that matches \eqref{eq:dinv_bijection} and \eqref{eq:inv_bijection} respectively when $l=0$ or $k=0$.

\begin{theorem}
    \label{thm:bijection}
    There is a bijection \[ \phi \colon \SW(\mu,k,l) \rightarrow \LD_0(\mu)^{\ast k, \bullet l} \] extending both \eqref{eq:dinv_bijection} and \eqref{eq:inv_bijection}.
\end{theorem}

\begin{proof}
    Let $\mu \vDash_0 n$. The bijection $\phi$ will be defined by induction on $m$, the maximal index in $\mu$ such that $\mu_m > 0$. Define $j=\mu_m$, and $\mu^- = (\mu_1, \dots, \mu_{m-1},0,\dots)$. As base case, we define $\phi(\varnothing) = \varnothing$, that is, the one segmented Smirnov word of size $0$ goes to the one path of size $0$.

    Recall that paths in $\LD_0(n)^{\ast k, \bullet l}$ can be easily characterized as being concatenation of paths of the form $N^iE^i$, in which all rises are decorated (see~\Cref{fig:path-example}, right). This precisely ensures that the area is zero, see~\Cref{def:area}.

    For a path in $\LD_0(\mu)^{\ast k, \bullet l}$, the maximal label $m$ occurs $j$ times by definition. Since $m$ is maximal, it can only occur as the top label of consecutive vertical steps. There are four distinct possibilities, illustrated in \Cref{fig:four_cases}, top:
    \begin{enumerate}
        \item $m$ labels a rise followed by a decorated valley.
        \item $m$ labels a rise followed by a non decorated valley (or is the last north step of the path).
        \item $m$ labels a north step on the diagonal, which is a decorated valley
        \item $m$ labels a north step on the diagonal, which is not decorated.
    \end{enumerate}
    To delete $m$ and get a path in $\LD_0(\mu^-)$, delete the north step carrying $m$ and the following east step. In the first case, also delete the rise decoration on that step and the valley decoration on the following north step; in the second and third case, also delete (rise or valley) decoration on that step. See~\Cref{fig:four_cases} for an illustration.

    \begin{figure}[!ht]
        \centering
        \includegraphics[width=0.8\textwidth]{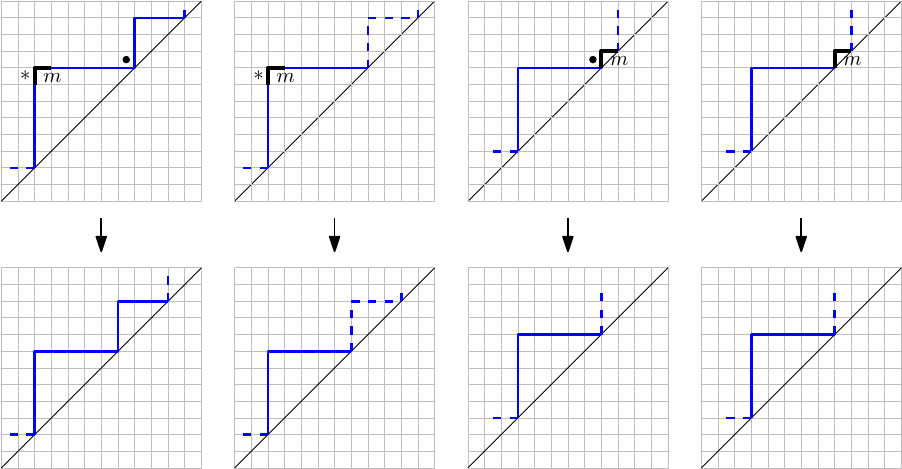}
        \caption{Four cases for the occurrence of the maximal label $m$, and how to delete it.
            \label{fig:four_cases}}
    \end{figure}

    We want to understand how to reverse this procedure, and insert the labels $m$ in a path in $\LD_0(\mu^-)$. Let us call \emph{blocks} of a decorated path with area equal to $0$ the subpaths separated by the non decorated north steps on the main diagonal: then all north steps of a block are decorated except the first one, so that a path in $\LD_0(n,k,l)$ will have $n-k-l$ blocks.

    By the inductive hypothesis, $\phi \colon \SW(\mu^-) \rightarrow \LD_0(\mu^-)$ is a bijection; we strengthen the inductive hypothesis by also asking that blocks of a Dyck path $D' \in \LD_0(\mu^-)$ correspond, from left to right, to blocks of $w' = \phi^{-1}(D') \in \SW(\mu^-)$, from right to left.

    We want to insert $j$ north steps labelled $m$ to get a path in $\LD_0(\mu)$. We do these insertions successively, according to the four types of such labels listed above.

    First, for each block of $D'=D^{(0)}$ that is not the last one, consider the last north step: one can insert a decorated rise labeled $m$, followed by an east step, and decorate the following valley: it was not decorated since it is between two blocks, and is contractible since it is preceded by at least two east steps. In total let $i$ be the number of insertions of this form; call the resulting path $D^{(1)}$. Each of these insertions joins two consecutive blocks of $D'$, and this corresponds to a \emph{peak insertion} in $w' = w^{(0)}$: call $w^{(1)}$ the word obtained by joining the corresponding blocks of $w'$, by replacing their block separators with an occurrence of $m$.

    Then, for each block of $D^{(1)}$, consider the last north step: one can insert a decorated rise labeled $m$, followed by an east step. Let $r-i$ be the number of such insertions, and $D^{(2)}$ be the resulting path. This corresponds to a \emph{double rise insertion} in the corresponding segmented Smirnov word, that is, inserting an occurrence of $m$ at the end of the corresponding blocks of $w^{(1)}$; call $w^{(2)}$ the word obtained this way. Notice that, if $l=0$, this corresponds to the bijection shown in \Cref{fig:op-bijection}.

    Next, for each block of $D^{(2)}$, one can extend it by adding a north step labelled $m$ on the diagonal, followed by an east step, and decorate this new north step which must be a contractible valley. Indeed, if the new valley is preceded by a rise, then it is always contractible; if not, since $m$ is maximal and there are no occurrences of $m$ on the main diagonal, then the new valley is necessarily contractible. Let $a-i$ be the number of insertions of this form, and call the resulting path $D^{(3)}$. This corresponds to a \emph{double fall insertion} in the corresponding segmented Smirnov word, that is, inserting an occurrence of $m$ at the beginning of the corresponding blocks of $w^{(2)}$; call $w^{(3)}$ the word obtained this way. Once again, if $k=0$, this corresponds to the bijection shown in \Cref{fig:op-bijection}.

    Finally, one can insert singleton blocks formed by a north step labeled $m$ on the diagonal followed by an east step among all blocks of $D^{(3)}$. Several such singleton blocks can be inserted between two consecutive blocks of $D^{(3)}$, as well as at the beginning and end of it. There must be $j-r-a+i$ such insertions, so that $j$ labels $m$ occur in the final path $D^{(4)}=D$. This corresponds to a \emph{valley insertion} in the corresponding segmented Smirnov word, that is, inserting an occurrence of $m$ as a singleton in between the corresponding blocks of $w^{(3)}$; call $w^{(4)} = w$ the word obtained this way. When $k=0$ or $l=0$, this also corresponds to the bijection shown in \Cref{fig:op-bijection}.

    We define $\phi(w) \coloneqq D$. Since each step of this construction is reversible, and we already showed in \Cref{thm:main} that segmented Smirnov words can be constructed this way (recall that double rise and double fall insertions commute), the thesis follows. Moreover, we showed that, when $k=0$ or $l=0$, this bijection restrict to the ones shown in \Cref{fig:op-bijection}, as desired.
\end{proof}

\subsection{Another statistic on segmented Smirnov words}
\label{sub:sdinv}

We now need to introduce a slightly tweaked version of our statistic on segmented Smirnov words, which better matches the combinatorics of the Delta conjecture. In order to do so, we need some preliminary definitions.

\begin{definition}
    Let $w \in \SW(n)$ and let $m \in \mathbb{N}$. For $1 \leq i \leq n$, let \[ \start_m(i) \coloneqq \max \{j \leq i \mid w_{j-1} > m \textrm{ or } j \text{ is initial} \}. \] We define \[ \height_m(i) \coloneqq \# \{ \start_m(i) \leq j < i \mid w_j < m \}, \] that is, the number of letters of $w$ strictly smaller than $m$ that are in between $w_i$ and the first thing to the left of $w_i$ that is either a letter strictly greater than $m$ or a block separator.
\end{definition}

\begin{definition}
    \label{def:smirnov_dinv}
    For a segmented Smirnov word $w$, we say that $(i,j)$ with $1 \leq i, j \leq n$ is a \emph{diagonal inversion} if $w_i > w_j$ and one of the following holds:
    \begin{enumerate}
        \item \label{dinv:right} $i$ is not a peak, $i < j$, $\height_{w_i}(i) = \height_{w_i}(j)$, and if $j=i+1$ then $j$ must be initial;
        \item \label{dinv:left} $i$ is not a peak, $i > j+1$, and $\height_{w_i}(i) = \height_{w_i}(j)+1$;
        \item \label{dinv:peak} $i$ is a peak and $(i,j)$ is a sminversion.
    \end{enumerate}
    We denote $\sdinv(w)$ the number of diagonal inversions of $w$, and \[ \overline{\SW}_{q}(\mu,k,l) \coloneqq \sum_{w \in \SW(\mu,k,l)} q^{\sdinv(w)}. \]
\end{definition}

\begin{remark}
    \label{rmk:dinv_is_inv}
    It is worth noticing that, if $i$ is not a double rise, then $(i,j)$ is a diagonal inversion if and only if it is a sminversion. Indeed, if $i$ is a valley or a double fall, then $i$ is either initial or a descent, and $\height_{w_i}(i) = 0$; it follows that, in \Cref{def:smirnov_dinv} condition~\ref{dinv:left} is never satisfied, and condition \ref{dinv:right}, namely $\height_{w_i}(j) = 0$, reduces to one of the four conditions in \Cref{def:smirnov_inv}. If $i$ is a peak instead, then diagonal inversions are sminversions by definition.
\end{remark}

\begin{example}
    Consider the word $231|3212|12$, as in \Cref{ex:sminv}. The only index that is a peak is $2$. For this index, we count sminversions to its right: $(2,5)$ and $(2,7)$. Let us compute the heights at all levels:
    \begin{center}
        \begin{tabular}{cl}
            the word    & $231|3212|12$                      \\
            the indices & $123\phantom{|}4567\phantom{|}89$  \\
            $\height_1$ & $000\phantom{|}0000\phantom{|}00$  \\
            $\height_2$ & $000\phantom{|}0001\phantom{|}01$  \\
            $\height_3$ & $011\phantom{|}0012\phantom{|}01.$
        \end{tabular}
    \end{center}
    So the pairs $(i,j)$ that form diagonal inversions are $(1,3)$, $(1,6)$, $(1,8)$, $(5,8)$, $(7,3)$, $(9,3)$, $(9,6)$ (for which $w_i=2$) and $(4,8)$ (for which $w_i=3$). Thus the $\sdinv$ of this word is $10$.
\end{example}
\begin{example}Let us compute $\sminv$ and $\dinv$ for all the elements of $\SW((2,1))$:
    \begin{center}
        \begin{tabular}{ccc}
            The word & its $\sminv$ & its $\dinv$ \\
            $1|1|2$  & 0            & 0           \\
            $1|2|1$  & 1            & 1           \\
            $2|1|1$  & 2            & 2           \\
            $1|12$   & 0            & 1           \\
            $1|21$   & 0            & 0           \\
            $12|1$   & 1            & 0           \\
            $21|1$   & 1            & 1           \\
            $121$    & 0            & 0
        \end{tabular}
    \end{center}
    So we obtain:
    \begin{align*}
         & \SW_q((2,1),0,0) = \overline{\SW}_q((2,1),0,0) = 1 + q + q^2 &  & \SW_q((2,1),1,0) = \overline{\SW}_q((2,1),1,0) = 1 + q  \\
         & \SW_q((2,1),1,1)=\overline{\SW}_q((2,1),1,1) = 1             &  & \SW_q((2,1),0,1) = \overline{\SW_q}((2,1),0,1) = 1 + q.
    \end{align*}
\end{example}

Indeed, by construction this statistic extends both $\inv$ and $\dinv$ on ordered multiset partitions; recall the map $\Pi$ defined before~\Cref{prop:sminv_to_inv}.

\begin{proposition}
    \label{prop:sdinv_to_dinv}
    The map $\Pi$ restricts to bijections \[ \SW(\mu, k, 0) \leftrightarrow \OP(\mu, n-k), \qquad \SW(\mu, 0, l) \leftrightarrow \OP(\mu, n-l) \] sending $\sdinv$ to $\dinv$ and $\inv$ respectively.
\end{proposition}

\begin{proof}
    \Cref{prop:sminv_to_inv} show that these are bijections, and since the statistics $\sminv$ and $\sdinv$ coincide when $k=0$, we only have to deal with the case $l=0$.

    Let $w \in \SW(\mu, k, 0)$. Since there are no descents, the blocks of $w$ are strictly increasing. It follows that, for every $i$, if $m \geq w_i$, then $\height_m(i) = h$ is the same as saying that $w_i$ is the $h\th$ smallest element in its block, and so the conditions in \Cref{def:smirnov_dinv} are the same as the ones in the definition of $\dinv$ on ordered multiset partitions. The thesis follows.
\end{proof}

\subsection{Equidistribution of \texorpdfstring{$\sminv$}{sminv} and \texorpdfstring{$\sdinv$}{sdinv}}
\label{sub:sminv_and_sdinv}

The statistics $\sminv$ and $\sdinv$ on segmented Smirnov words are equidistributed. Indeed, if we let \[ \overline{\SW}_{x;q}(n,k,l) \coloneqq \sum_{w \in \SW(n,k,l)} q^{\sdinv(w)} x_w, \] we have the following.

\begin{theorem}
    \label{thm:sdinv_main}
    For any $n,k,l$ with $k+l < n$, we have the identity \[ \SW_{x;q}(n,k,l) = \overline{\SW}_{x;q}(n,k,l). \]
\end{theorem}

As before, we split the proof in several lemmas.

\begin{lemma}
    \label{lem:sdinv_preserve_inversions}
    Let $w \in \SW(n)$, and let $m \geq \max w$. Let $\overline{w} \in \SW(n+1)$ be a word obtained from $w$ by inserting an occurrence of $m$ as a double fall, double rise, or singleton. Then, if $(i,j)$ forms a diagonal inversion in $w$, the corresponding indices in $\overline{w}$ also form a diagonal inversion. If instead $m$ is inserted as a peak, the same holds as long as $m > \max w$.
\end{lemma}

\begin{proof}
    Suppose that $(i,j)$ is a diagonal inversion in $w$. If $m$ is not inserted as a peak, then $\height_a$ does not change on any letter of $w$ for all $a \leq m$, so if two letters of $w$ form a diagonal inversion, they will also form a diagonal inversion in $\overline{w}$. If $m$ is inserted as a peak instead (i.e.\ replacing a block separator), then $\height_a$ does not change on any letter of $w$ for all $a < m$, but since every letter of $w$ is strictly smaller than $m$, then if two letters of $w$ form a diagonal inversion, they will also form a diagonal inversion in $\overline{w}$.
\end{proof}

By \Cref{rmk:dinv_is_inv} peaks, double falls, and valleys, sminversions and diagonal inversions coincide. It follows that \Cref{lem:peak_insertion}, \Cref{lem:fall_insertion}, and \Cref{lem:valley_insertion} also hold for $\sdinv$, the proof being exactly the same. For the double rise insertion, we will need a different strategy. The reader is invited to follow along using \Cref{ex:sdinv-insertion}.

\begin{lemma}
    \label{lem:rise_sdinv_insertion}
    Let $w \in \SW(n,k,l)$, and let $m \geq \max w$ with no final $m$ in $w$. The $q$-enumerator with respect to $\sdinv$ of words obtained from $w$ by inserting $s$ occurrences of $m$ in $w$ as a double rise (i.e.\ as the final element of a block) is \[ q^{\binom{s}{2}} \qbinom{n-k-l}{s} q^{\sdinv(w)}. \]
\end{lemma}

\begin{proof}
    We proceed in the same fashion as in \cite{Wilson2016Equidistribution}*{Subsection~4.3}. Let $w = (w^1, \dots, w^h)$, and let us preliminarily define $\lvert w^i \rvert_m \coloneqq \# \{ a \mid w^i_a < m \}$, that is, the number of letters strictly smaller than $m$ in $w^i$. We sort the blocks lexicographically by $(-\lvert w^i \rvert_m, i)$: we say that $w^j \preccurlyeq w^i$ if either $\lvert w^i \rvert_m < \lvert w^j \rvert_m$, or $\lvert w^i \rvert_m = \lvert w^j \rvert_m$ and $j < i$; in other words, if we disregard all the occurrences of $m$, the leftmost biggest block is the first in this order, and the rightmost smallest block is the last.

    We claim that adding an occurrence of $m$ at the end of $w^i$ creates $\# \{ j \mid w^j \prec w^i \}$ new diagonal inversions. Notice that, by definition, for each $i$ the values of $\height_m$ on the indices (in $w$) corresponding to elements in $\{ a \mid w^i_a < m \}$ are exactly the interval $[0, \lvert w^i \rvert_m - 1]$, as $w^i$ has no entry strictly greater than $a$. Also notice that, if we insert an occurrence of $m$ at the end of $w^i$, its $\height_m$ is going to be exactly $\lvert w^i \rvert_m$.

    Now, insert an occurrence of $m$ at the end of $w^i$. Take $j$ such that $w^j \prec w^i$. If $\lvert w^i \rvert_m < \lvert w^j \rvert_m$, then in $w^j$ there are two letters, both strictly smaller than $m$, whose indices have $\height_m$ equal to $\lvert w^i \rvert_m$ and $\lvert w^i \rvert_m - 1$ respectively, so whether $i < j$ or $i > j$, said occurrence of $m$ forms a diagonal inversion with some entry in $w^j$. If $\lvert w^i \rvert_m = \lvert w^j \rvert_m$, then by definition $i > j$ and the index of the last letter of $j$ has $\height_m$ equal to $\lvert w^i \rvert_m - 1$, so it forms a diagonal inversion with the occurrence of $m$ at the end of $i$. It is easy to see that no elements in $w^j$ with $w^i\preccurlyeq w^j$ create diagonal inversions with the newly inserted letter $m$: they all have $\height_m$ strictly smaller that $|w^i|-1$, except for the letter immediately preceding the newly inserted $m$, with which does not create a diagonal inversion by definition.

    Summarizing, said ordering of the blocks of $w$ is such that inserting an occurrence of $m$ at the end of the $a^\mathsf{th}$ block creates $a-1$ new diagonal inversions; since we have to insert $s$ such occurrences, and we have to do that in distinct blocks, the new diagonal inversions are $q$-counted by $q^{\binom{s}{2}} \qbinom{n-k-l}{s}$, which concludes the proof.
\end{proof}

\begin{example}
    \label{ex:sdinv-insertion}
    Let us consider a segmented Smirnov word with maximal letter $3$, such that $3$ is never a final element. We record $\height_3$ for every letter, $\lvert w^i \rvert_3$ for every block and indicate the ordering on the blocks induced by $\prec$:
    \begin{center}
        \begin{tabular}{rl}
            the word    & $21\;|\;121\;|\;321\;|\;2\;|\;2132\;|\;12$                                                    \\
            $\height_3$ & $01\;\phantom{|}\; 012\;\phantom{|}\;001\;\phantom{|}\;0\;\phantom{|}\;0122\;\phantom{|}\;01$ \\
            $\lvert w^i \rvert_3$   &
            \framebox[\widthof{21\;}]{2}\phantom{|}\framebox[\widthof{\;121\;}]{3}\phantom{|}\framebox[\widthof{\;321\;}]{2}\phantom{|}\framebox[\widthof{\;2\;}]{1}\phantom{|}\framebox[\widthof{\;2132\;}]{3}\phantom{|}\framebox[\widthof{\;12}]{2}
            \\[5pt]ordering of the blocks&
            \framebox[\widthof{21\;}]{2}\phantom{|}\framebox[\widthof{\;121\;}]{0}\phantom{|}\framebox[\widthof{\;321\;}]{3}\phantom{|}\framebox[\widthof{\;2\;}]{5}\phantom{|}\framebox[\widthof{\;2132\;}]{1}\phantom{|}\framebox[\widthof{\;12}]{4}
        \end{tabular}.
    \end{center}
    Let us insert a $3$ as a final element into block number $3$, with respect to $\prec$ (the letter in red):
    \begin{center}
        \begin{tabular}{rl}
            the word    & $2\textcolor{blue}{1}\;|\; 1\textcolor{blue}{2}1\;|\;321\textcolor{red}{3}\;|\;2\;|\;213\textcolor{blue}{2}\;|\;12$                                                   \\
            $\height_3$ & $0\textcolor{blue}{1}\;\phantom{|}\;0\textcolor{blue}{1}2\;\phantom{|}\;001\textcolor{red}{2}\;\phantom{|}\;0\;\phantom{|}\;012\textcolor{blue}{2}\;\phantom{|}\;01$.
        \end{tabular}
    \end{center}
    This newly inserted letter creates diagonal inversions only with $3$ letters: exactly one in each of the preceding (with respect to $\prec$) blocks (indicated in blue).
\end{example}

\Cref{thm:sdinv_main} now follows easily.

\begin{proof}[Proof of \Cref{thm:sdinv_main}]
    We prove the equivalent statement \[ \langle \SF(n,k,l), h_\mu \rangle = \overline{\SW}_q(\mu, k, l). \] By the same argument used in the proof of \Cref{thm:main}, using \Cref{lem:peak_insertion}, \Cref{lem:rise_sdinv_insertion}, \Cref{lem:fall_insertion}, and \Cref{lem:valley_insertion}, it follows that $\overline{\SW}_q(\mu, k, l)$ satisfies the same recursion as $\SW_q(\mu, k, l)$ with the same initial conditions, so the two must coincide.
\end{proof}

\begin{remark}
    \label{rmk:peak_insertion}
    \Cref{prop:sdinv_to_dinv} holds regardless of the ordering on the blocks we use for the peak insertion: indeed, when $k=0$ or $l=0$, peak insertions never happen, so this adds a layer of freedom to the definition of $\sdinv$ while maintaining the same essential properties. We decided to keep the ordering linear for simplicity, but since every possible ordering defines a statistic with the same distribution; it is possible that a different ordering yields a variant of this statistic that has a simpler or more uniform description.
\end{remark}

\subsection{Unified Delta theorem at \texorpdfstring{$t=0$}{t=0}}
\label{sub:unified_delta_theorem}

Since our $\sdinv$ statistic reduces to the $\dinv$ and $\inv$ statistics on ordered multiset partitions when either $l=0$ or $k=0$, the bijection in \Cref{thm:bijection} extends \eqref{eq:dinv_bijection} and \eqref{eq:inv_bijection} also with respect to the statistic. It follows that, when $k=0$ or $l=0$, our bijection maps the $\sdinv$ of a segmented Smirnov word to the $\dinv$ of the corresponding decorated Dyck path; therefore, this bijection produces a $q$-statistic on labelled Dyck paths with both decorated rises and decorated contractible valleys that matches the expected symmetric function expression when $t=0$, establishing a very solid first step towards the statement of a unified Delta conjecture.

In particular, we know that \[ \left. \Theta_{e_k} \nabla e_{n-k} \right\rvert_{t=0} = \sum_{\pi \in \OP(n, n-k)} q^{\dinv(\pi)} x_\pi = \sum_{D \in \LD_0(n)^{\ast k}} q^{\dinv(D)} x^D \] and that \[ \left. \Theta_{e_l} \nabla e_{n-l} \right\rvert_{t=0} = \sum_{\pi \in \OP(n, n-l)} q^{\inv(\pi)} x_\pi = \sum_{D \in \LD_0(n)^{\bullet l}} q^{\dinv(D)} x^D, \] so, now that we have a unified combinatorial model for $\left. \Theta_{e_k} \Theta_{e_l} \nabla e_{n-k-l} \right\rvert_{t=0}$, it is natural to try to interpret it in the framework of the Delta conjecture. We achieve that with the following definition.

\begin{definition}
    For $D \in\LD_0(n)^{\ast k, \bullet l}$, we define $\dinv(D) \coloneqq \sdinv(\phi^{-1}(D))$.
\end{definition}

While not entirely explicit, this statistic solves the long-standing problem of giving a unique model for the two versions of the Delta conjecture, at least when $t=0$; indeed, \Cref{thm:unified_Delta_t_equals_0} now follows as a trivial corollary of \Cref{thm:bijection} and \Cref{thm:sdinv_main}. It would be very interesting to have an explicit description of the statistic on Dyck paths, maybe using a different variant of the peak insertion (as per \Cref{rmk:peak_insertion}), in order to possibly derive an extension to objects with positive area.

\section{`Catalan' case via fundamental quasisymmetric expansion}
\label{sec:fundamental-quasisym}

In this section, we give an expansion for $\SW_{x;q}(n,k,l)$ in terms of fundamental quasisymmetric functions. This is a compact way to rewrite the expansion into segmented Smirnov words. We then use it to compute $\langle \SF(n,k,l), e_n \rangle$, which is the coefficient of $s_{1^n}$ in the Schur expansion of $\SF(n,k,l)$. Is it thus the (conjectured) graded dimension of the sign-isotypical component of the $(1,2)$-coinvariant module. This is often referred to as the `Catalan' case, as, for the classical shuffle theorem, $\langle \nabla e_n, e_n \rangle$ yields a $(q,t)$-analogue of the Catalan numbers.

\subsection{Fundamental quasisymmetric expansion}
\label{sub:fundamental}
We need some definitions before stating our expansion. 

\begin{definition}
    Let $w$ be a segmented Smirnov word. For $1 \leq i \leq n$, we say that $i$ is \emph{thick} if $i$ is initial, or $i$ is not initial and $w_{i-1} > w_i$ (i.e.\ $i-1$ is a descent); we say it is \emph{thin} otherwise, that is, $i$ is not initial and $w_{i-1} > w_i$ (i.e.\ $i-1$ is an ascent).
\end{definition}

\begin{definition}
    Let $\sigma$ be a segmented permutation of size $n$, $i<n$, and let $j$ be such that $\sigma_j = \sigma_i + 1$ (so $j=\sigma^{-1}(\sigma_i+1)$). We say that $\sigma_i$ is \emph{splitting} for $\sigma$ if either of the following holds:
    \begin{itemize}
        \item $i$ is thick and $j$ is thin;
        \item $i$ and $j$ are both thin and $i < j$;
        \item $i$ and $j$ are both thick and $j < i$.
    \end{itemize}
\end{definition}

Let $\Split(\sigma) = \{ 1 \leq i \leq n-1 \mid i \text{ is splitting for } \sigma \}$. For example, if $\sigma=152|93|6487$, then $\Split(\sigma) = \{4, 7\}$. Recall that, for any subset $S \subseteq [n-1]$, the fundamental quasisymmetric function $Q_{S,n}$ is defined as
\begin{equation}
    \label{eq:fundamental}
    Q_{S,n}=\sum_{\substack{1 \leq v_1 \leq v_2 \leq \ldots \leq v_n \\ v_j < v_{j+1} \text{ if } j \in S}} x_{v_1} x_{v_2} \cdots x_{v_n}.
\end{equation}
The main result of this section is the following.

\begin{proposition}
    \label{prop:fundamental}
    \[ \SW_{x;q}(n,k,l) = \sum_{\sigma \in \SW(1^n,k,l)} q^{\sminv(\sigma)} Q_{\Split(\sigma), n}. \]
\end{proposition}

Define the \emph{reading order} of $w$ by scanning the entries from the smallest value to the biggest, and for a given value we first read the thin entries from right to left, and then the thick ones from left to right. Explicitly, it is the linear order $\leq_w$ on $[n]$ defined by $i <_w j$ if either of the following holds:
    \begin{enumerate}
        \item $w_i < w_j$;
        \item $w_i = w_j$, $i$ is thin and $j$ is thick;
        \item $w_i = w_j$, $i$ and $j$ are both thin and $i > j$;
        \item $w_i = w_j$, $i$ and $j$ are both thick and $i < j$.
    \end{enumerate}
    The \emph{standardization} $\st(w)$ of $w$ is then the segmented permutation $\sigma$ of the same shape such that $\leq_w = \leq_\sigma$, that is $\sigma^{-1}(1) <_w \sigma^{-1}(2) <_w \ldots <_w \sigma^{-1}(n)$.

\begin{example}
    Let $w = {1}\ul{2}{1}|{31}|{21}\ul{3}{2}$ where thin entries are underlined. Then $\st(w) = 152|93|6487$.
\end{example}

\begin{lemma}
    For any segmented Smirnov word $w$, we have that $w$ and $\st(w)$ have the same ascents, descents, and sminversions.
\end{lemma}

\begin{proof}
    Let us write $\sigma=\st(w)$. The fact that $w$ and $\sigma$ have the same ascents and descents follows immediately from the fact that $w_i < w_j \implies i <_w j$. 

    We now want to show that $(i,j)$ is a sminversion in $w$ if and only if it is a sminversion in $\sigma$. Notice first that $i$ is thick in $w$ if and only if it is thick in $\sigma$, and that for $i < j$ and $w_i = w_j$, we have $\sigma_i < \sigma_j$ if and only if $j$ is thick.

    Suppose that $(i,j)$ is a sminversion in $w$: then $w_i > w_j$, so by construction $\sigma_i > \sigma_j$. Now if $j$ is initial in $w$, then it is in $\sigma$ as well. If $w_{j-1} > w_i$, then that is the case for $\sigma$ as well. This deals with the cases \eqref{inv:first_letter} and \eqref{inv:bigger_entry} in \Cref{def:smirnov_inv}. Now cases \eqref{inv:equal_then_first_letter} and \eqref{inv:equal_then_bigger_entry} can be rewritten as $i < j-1$, $w_{j-1} = w_i$, and $j-1$ is thick. This implies that $j-1$ will come after $i$ in the reading order, and thus $\sigma_{j-1} > \sigma_i$. In each case $(i,j)$ is a sminversion in $\sigma$, as desired.

    Suppose now that $(i,j)$ is a sminversion in $\sigma$. This means that $\sigma_i > \sigma_j$ and either $j$ is initial or $\sigma_{j-1} > \sigma_i > \sigma_j$. In either case $j$ is thick, and in particular we must have $w_i > w_j$. If $j$ is initial in $\sigma$, then it is in $w$ as well. And if $\sigma_{j-1} > \sigma_i$, then either $w_{j-1} > w_i$, or $w_{j-1} = w_i$ and $j-1$ is thick. In either case $(i,j)$ is a sminversion in $w$, as desired.
\end{proof}
%
%
%
%

\begin{proof}[Proof of~\Cref{prop:fundamental}]
    Thanks to the previous lemma, it is enough to show that, for any $\sigma \in \SW(1^n,k,l)$, we have \[ Q_{\Split(\sigma), n} = \sum_{\substack{w \in \SW(n,k,l) \\ \st(w)=\sigma}} x_w. \]
    We fix $\sigma$. For $w$ a segmented Smirnov word, write $x_w = x_{v_1} x_{v_2} \cdots x_{v_n}$ where $v_j$ is the increasing rearrangement of $w_i$'s according to the linear order $\leq_w$. Comparing the formula above with~\eqref{eq:fundamental}, we are led to show the following claim: $w$ is a segmented Smirnov word such that $\st(w) = \sigma$ if and only if $w$ is a segmented word with the same shape as $\sigma$, and for each $i, j$ such that $\sigma_j = \sigma_i + 1$, we have $w_i \leq w_j$, and $w_i < w_j$ if $i \in \Split(\sigma)$. 

    Suppose first that $w \in \SW(n)$ satisfies $\st(w) = \sigma$. By definition of standardization $w$ and $\sigma$ have the same shape. If $\sigma_j = \sigma_i + 1$, in particular $\sigma_i < \sigma_j$, so necessarily $w_i \leq w_j$. If moreover $w_i = w_j$, then by definition $i$ is not splitting in $\sigma$, and thus $i \in \Split(\sigma)$ implies $w_i < w_j$.

    Conversely, suppose that $w$ is a segmented word with the same shape as $\sigma$, such that if $\sigma_j = \sigma_i + 1$, then $w_i \leq w_j$, and $w_i < w_j$ if $i \in \Split(\sigma)$. This implies the Smirnov condition for $w$: indeed, assume we have $w_i=w_{i+1}$ in a block for some $i$. In the order $\leq_\sigma$ consider the chain between $i$ and $i+1$. No element but the last one can be split by the hypothesis, so the chain consists of thin elements with decreasing values followed by thick elements with increasing values. But no two elements of such a chain can be adjacent in a block, as is checked by direct inspection. This contradicts the hypothesis, and thus $w\in\SW(n)$.
    
    We must show $\st(w)=\sigma$, that is, the reading order $\leq_w$ coincides with $\leq_\sigma$. Assume then $i<_\sigma j$. This first entails $w_i\leq w_j$ by hypothesis. If $w_i<w_j$ then $i<_w j$ by definition of $<_w$. Assume now $w_i=w_j$; we can suppose $i$ and $j$ consecutive in $<_\sigma$, that is, $\sigma_j = \sigma_i + 1$. By the hypothesis we have $i\notin \Split(\sigma)$: since thin or thick indices coincide in $w$ and $\sigma$ by the hypothesis, the definition of the reading order implies that $i<_w j$, as desired.
\end{proof}


\subsection{The `Catalan' case}

In this section, we look at the combinatorics of $\langle \Theta_{e_k} \Theta_{e_l} \Ht_{(n-k-l)}, e_n \rangle \rvert_{t=0}$. As it is happens often, this case is significantly simpler, and we are able to provide an explicit formula in terms of $q$-binomials.

\begin{theorem}
    \label{thm:catalan}
    \[ \langle \SF(n,k,l), e_n \rangle = q^{\binom{n-k-l}{2}} \qbinom{n-1}{k+l} \qbinom{k+l}{k}. \]
\end{theorem}

The product of $q$-binomial coefficients above clearly enumerates sequences $v$ of length $n-1$ on the alphabet $\{0,1,2\}$ that have $n-k-1$ occurrences of $0$, $k$ occurrences of $1$ and $l$ occurrences of $2$, according to inversions $i<j$ such that $v_i > v_j$. The inversion statistic is called {\em mahonian} on $\{0,1,2\}$. Our starting point is the following.

\begin{proposition}
    \label{prop:catalan-sw}
    \[ \langle \SF(n,k,l), e_n \rangle = \sum_{\substack{\sigma \in \SW(1^n,k,l)\\ \Split(\sigma) = [n-1]}} q^{\sminv(\sigma)}.\] 
\end{proposition}

\begin{proof}
    Recall that, if $f$ is symmetric, then its fundamental expansion determines its Schur expansion as follows (see \cites{Stanley1999EnumerativeCombinatoricsV2, Gessel2019}) \[ f = \sum_{\alpha \vDash n} c_\alpha F_\alpha \implies \sum_{\alpha \vDash n} c_\alpha s_\alpha. \] 
    Here $F_\alpha = Q_{S(\alpha),n}$ where $S(\alpha)$ is the subset of $[n-1]$ given by the partial sums of $\alpha$, and $s_\alpha$ is given by the Jacobi-Trudi formula, and is thus equal to $0$ or a Schur function (up to sign). Now $Q_{[n-1],n} = F_{1^n}$ and the only composition $\alpha$ such that $s_\alpha = \pm s_{1^n}$ is $1^n$. Since $e_n = s_{1^n}$, and the Schur basis is orthonormal, $\langle \SF(n,k,l), e_n \rangle$ is the coefficient of $Q_{[n-1],n}$ in the expansion of~\Cref{prop:fundamental}, which is the desired expression.
\end{proof}

We now need to characterize all the segmented permutations $\sigma \in \SW(1^n,k,l)$ such that $\Split(\sigma) = [n-1]$. Let us use the notation $B=n-k-l$ for the number of blocks.
 
\begin{proposition} \label{prop:characterization-allsplit}
    A segmented permutation $\sigma\in \SW(1^n,k,l)$ satisfies $\Split(\sigma) = [n-1]$ if and only if :
    \begin{itemize}
        \item $\sigma_1 = B+l$;
        \item the values $1,2,\ldots,B+l$ occur in $\sigma$ from right to left, and are in thick positions, with $B$ of them initial;
        \item the values $B+l+1,\ldots,n$ occur in $\sigma$ from left to right in thin positions.
    \end{itemize}
\end{proposition}
\begin{example}
    The segmented permutation $\sigma=67|584|3|291$ is an element of $\SW(1^9,3,2)$ such that $\Split(\sigma)=[8]$.
\end{example}
\begin{proof}[Proof of Proposition~\ref{prop:characterization-allsplit}]
    The condition is sufficient: if $\sigma$ is constructed as described, all positions are clearly splitting.

    Conversely, assume $\Split(\sigma) = [n-1]$. By definition, since every index $i$ is splitting, if $\sigma_j = \sigma_i + 1$, if $i$ is thin then $j$ must also be thin; it follows that we can't have a thick letter that is greater than a thin one, so the thick letters must be $1, 2, \dots, B+l$. Moreover, thick letters must appear in increasing order from right to left, and thin letters must appear in increasing order from left to right. We must have $\sigma_1=B+l$ in order for it to be thick, and of course $B$ positions must be initial since $\sigma$ has $B$ blocks. This gives the desired characterization.
\end{proof}

  For $w \in \SW(n)$, we define its \emph{skeleton} to be $\skel(w)$ by $\skel(w)_i = 0$ (resp.\ $1$, resp.\ $2$) if $i$ is an initial element (resp.\ $i$ is thick and not initial, resp.\ $i$ is thin). Clearly $\skel(w)$ starts with $0$, since $1$ is always an initial position. By removing the initial $0$ and changing $0,1,2$ to $\mid, >, <$ respectively, it is apparent that the skeleton encodes precisely the shape, ascents and descents of $w$. 

\begin{example}
    If we take again $\sigma=67|584|3|291$, then $\skel(\sigma)=020210021$, which  can be rewritten as \[ \cdot < \cdot \mid \cdot < \cdot > \cdot \mid \cdot \mid \cdot < \cdot > \cdot \] in this notation. Conversely, if we start with this skeleton, and insert the numbers from $1$ to $9$ in thick positions from right to left first, and then in thin positions from left to right, we get back $\sigma=67|584|3|291$.
\end{example}
  
\begin{corollary}
    \label{cor:characterization-allsplit}
    The map $\skel$ restricts to a bijection between $\{ \sigma \in \SW(1^n) \mid \Split(\sigma) = [n-1] \}$ and $\{ v \in \{0,1,2\}^n \mid v_1 = 0 \}$. 
\end{corollary}

Indeed, any choice of $v \in \{ v \in \{0,1,2\}^n \mid v_1 = 0 \}$ determines a unique such $\sigma$ thanks to \Cref{prop:characterization-allsplit}. This result implies that the (conjectured) dimension of the sign-isotypical component is $3^{n-1}$, as expected \cite{Bergeron2020}*{Table~3}. \Cref{thm:catalan} further refines this statement.

\begin{proof}[Proof of \Cref{thm:catalan}]
    Using \Cref{prop:catalan-sw}, we have to prove that (recall $B=n-k-l$)
    \begin{equation}
        \label{eq:mahonian_to_prove}
        \sum_{\substack{\sigma \in \SW(1^n,k,l) \\ \Split(\sigma) = [n-1]}} q^{\sminv(\sigma)} = q^{\binom{B}{2}} \qbinom{n-1}{k+l} \qbinom{k+l}{k}. 
    \end{equation}
    Let $\sigma \in \SW(1^n,k,l)$ such that $\Split(\sigma) = [n-1]$. In any sminversion $(i,j)$, $j$ is necessarily thick by~\Cref{def:smirnov_inv}. Let $j$ be a thick index in $\sigma$. If $j$ is initial, any $i<j$ forms a sminversion of $\sigma$ by~\Cref{prop:characterization-allsplit}. If not, we have $w_{j-1}>w_j$, and then $i<j$ forms a sminversion if and only if $i<j-1$ and $j-1$ is thin, using~\Cref{prop:characterization-allsplit} again.
    
    This translates to skeletons as follows: if $\skel(\sigma) = v$, then sminversions correspond to pairs $i<j$ in $[n]$ such that either $v_j=0$, or $i < j-1$ and $v_{j-1} v_j = 21$. Writing $v = 0v'$, and counting the $\binom{B}{2}$ pairs $(0,0)$ separately, we can rewrite this as \[\sminv(\sigma) = \binom{B}{2} + \inv_{(0,1),(0,2)}(v') + \maj_{(1,2)}(v'),\] where $\inv_{(0,1),(0,2)}(v')$ is the number of indices $i<j$ in $v'$ such that $v_i \in \{1,2\}$ and $v_j=0$, and $\maj_{(1,2)}(v')$ is the sum of all $i<n$ such that $v'_i=2$ and $v_{i+1}=1$.
    
    It turns out that $\inv_{(0,1),(0,2)}+\maj_{(1,2)}$ is an example of a mixed statistic studied in~\cite{Kasraoui2009}. In particular, \cite{Kasraoui2009}*{Corollary~1.10} with $X = \{0,1,2\}$, $g(0)=g(2)=\infty$ and $g(1)=2$ shows that the statistic is mahonian; see the remark after~\Cref{thm:catalan}. By~\Cref{cor:characterization-allsplit}, this gives us the desired product of binomial coefficients, and proves our result.    
\end{proof}
    
Kasraoui's result above is in fact proved bijectively. It is based on Foata's \emph{second fundamental transformation}; see for instance \cite{Lothaire1997}*{\S 10.6}. Let us describe this construction in our case: we will define a function $\gamma$ from words in $X$ to itself, by induction on length. $\gamma$ sends the empty word to itself, and $\gamma(a)=a$ for $a\in X$. Assume now $v$ has length $\geq 2$. If $v=v'0$, define $\gamma(v)=\gamma(v')0$; if $v=v'2$, define $\gamma(v)=\gamma(v')2$. Now assume $v=v'1$. If $\gamma(v')=v''2$, define $\gamma(v)=2v''1$. Finally, if $\gamma(v')=v''a$ with $a\in\{0,1\}$, factor $\gamma(v')=2^{k_1} b_1 2^{k_1} b_2\cdots 2^{k_j} b_j$ where $b_i\in\{0,1\}$ (with $b_j=a$) and $k_i\geq 0$ for all $i$. Then define $\gamma(v)=b_1 2^{k_1} b_2 2^{k_2}\cdots b_j 2^{k_j} 1.$
The claim is that $\gamma$ is a bijection such that for any $v$, the number of inversions of $\gamma(v)$ is $\inv_{(0,1),(0,2)}(v) + \maj_{(1,2)}(v)$. We leave the proof to the reader, the general case being done in~\cite{Kasraoui2009}.

\section{The maximal case \texorpdfstring{$k+l=n-1$}{k+l=n-1}}

We focus in this last section on various aspects of the case $k+l=n-1$ of \Cref{thm:main}. The combinatorial side now involves only Smirnov words. It is also conjecturally giving the graded Frobenius characteristic of the subspace of the $(1,2)$-coinvariant space of maximum total degree in the fermionic variables $\boldsymbol{\zeta}_n, \boldsymbol{\xi}_n$ (cf. \eqref{eq:thetatheta-module}), extending (at least conjecturally) the results in \cite{KimRhoades2022NoncrossingPartitions}.

It turns out that this special case is linked to familiar instances of geometric and combinatorial constructions.

\subsection{Chromatic quasisymmetric functions}
\label{sec:chromatic}

Smirnov words are in bijection with proper colorings of the path graph, which means that our symmetric function, suitably specialized, coincides with the chromatic quasisymmetric function of the path of length $n$, which also coincides with the Frobenius characteristic of the representation of $\mathfrak S_n$ on the cohomology of the toric variety $\mathcal{V}_n$ associated with the Coxeter complex of type $A_{n-1}$ (see \cite{ShareshianWachs2016Chromatic}*{Section~10}).

Given a graph $G = (V,E)$, a \emph{proper coloring} is a function $c \colon V \rightarrow \mathbb{Z}_+$ such that $\{i, j\} \in E \implies c(i) \neq c(j)$. If $V = [n]$, a \emph{descent} of a coloring is an edge $\{i,j\} \in E$ such that $i < j$ and $c(i) > c(j)$. The chromatic quasisymmetric function of $G$ is defined as \[ \chi_G(X;q) = \sum_{\substack{c \colon V \rightarrow \mathbb{Z}_+ \\ c \text{ proper}}} q^{\mathsf{des}(c)} \prod_{v \in V} x_{c(v)}, \] where $\mathsf{des}(c)$ is the number of descents of $c$.

For the path graph $G_n = 1-2-\dots-n$, proper colorings clearly correspond to Smirnov words of length $n$ via $c\mapsto c(1) c(2) \dots c(n)$. It follows from \Cref{thm:main} that \[ X_{G_n}(\mathbf{x};u) = \left. \sum_{k=0}^{n-1} u^k \Theta_{e_k} \Theta_{e_{n-k-1}} e_1 \right\rvert_{q=1, t=0}.\]

This suggests the existence of an extra $q$-grading on the cohomology of the permutahedral toric variety $\mathcal{V}_n$: indeed the graded Frobenius characteristic of that cohomology is known to be given by $\omega X_{G_n}(\mathbf{x};u)$, see~\cite{ShareshianWachs2012ChromaticHessenberg}.

\subsection{Parallelogram polyominoes}
\label{sec:polyominoes}

\begin{figure}
    \begin{minipage}{.5\textwidth}
        \centering
        \includegraphics{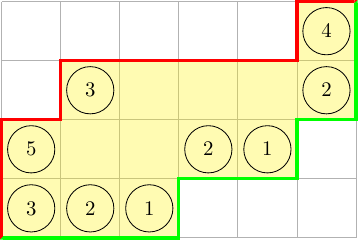}
    \end{minipage}%
    \begin{minipage}{.5\textwidth}
        \centering
        \includegraphics{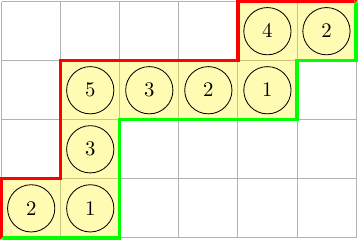}
    \end{minipage}
    \caption{Two labelled parallelogram polyominoes of size $6 \times 4$, the right one has area $0$.}
    \label{fig:polyomino}
\end{figure}

A \emph{parallelogram polyomino} of size $m \times n$ is a pair of north-east lattice paths on a $m \times n$ grid, such that first is always strictly above the second, except on the endpoints $(0,0)$ and $(m,n)$. A labelling of a parallelogram polyomino is an assignment of positive integer labels to each cell that has a north step of the first path as its left border, or an east step of the second path as its bottom border, such that columns are strictly increasing bottom to top and rows are strictly decreasing left to right. See \Cref{fig:polyomino} for two examples of labelled parallelogram polyominoes. In \cite{DAdderioIraciLeBorgneRomeroVandenWyngaerd2022TieredTrees} it is conjectured that $\Theta_{e_{m-1}} \Theta_{e_{n-1}} e_1$ enumerates labelled parallelogram polyominoes of size $m \times n$ with respect to two statistics, one of which is (a labelled version of) the area, and the other is unknown.

It is immediate to see that parallelogram polyominoes of size $(n-k) \times (k+1)$ and area $0$ are again in bijection with Smirnov words of length $n$ with $k$ ascents. Indeed, reading the labels of such a polyomino bottom to top, left to right, yields a Smirnov word of size $n$ with $k$ ascents, and the correspondence is bijective. For example, the Smirnov word corresponding to the area $0$ polyomino on the right in \Cref{fig:polyomino} is $213532142$.

In particular, sminversions on Smirnov words define a statistic on this subfamily of parallelogram polyominoes, proving the conjectural identity and partially answering Problem~7.13 from \cite{DAdderioIraciLeBorgneRomeroVandenWyngaerd2022TieredTrees} in the case when the area is $0$.

\subsection{The case \texorpdfstring{$q=0$}{q=0}}
\label{sec:topdegree}

Note that in this case, it is known~\cite{IraciRhoadesRomero2023FermionicTheta} that the symmetric function in \Cref{thm:main} is the Frobenius characteristic of the $(0,2)$-case. It was also shown that the high-degree part of this module has a basis indexed by noncrossing partitions \cite{KimRhoades2022NoncrossingPartitions}. In particular, this means that there is a bijection between segmented permutations with one block (that is, permutations) with zero $\sminv$, and noncrossing partitions. It would be interesting to see if, using this combinatorial model, one could get a combinatorial basis for the $(1,2)$-coinvariant ring.

\begin{lemma}
    \label{lem:231}
    Permutations with zero $\sminv$ are exactly $231$-avoiding permutations, that is, permutations $\sigma$ with no $i<j<k$ such that $\sigma_k<\sigma_i<\sigma_j$.
\end{lemma}

\begin{proof}
    Let $\sigma$ be a permutation, and suppose that it has a $231$ pattern, that is, that there exist indices $i < j < k$ such that $\sigma_k < \sigma_i < \sigma_j$. Let $m = \min \{j < a \leq k \mid \sigma_m < \sigma_i\}$; by definition, $i < j \leq m-1$, and $\sigma_{m-1} > \sigma_i$, so $(i,m)$ is a sminversion of $\sigma$. It follows that permutations with zero $\sminv$ are $231$-avoiding permutations. Since a sminversion in a permutation corresponds to a $231$ pattern, this concludes the proof.
\end{proof}

Let $\pi$ be a noncrossing partition, and let $\phi(\pi)$ be the permutation that, in one line notation, is written by listing the blocks of $\pi$ sorted by their smallest element, with the elements of each block sorted in decreasing order. Let us call \emph{decreasing run} of a permutation $\sigma$ a maximal subsequence of consecutive decreasing entries of $\sigma$ (in one line notation): then the blocks of $\pi$ correspond to the decreasing runs of $\phi(\pi)$. For instance, if $\pi = \{ \{1,2,5\}, \{3,4\}, \{6,8,9\}, \{7\}\}$, then $\phi(\pi) = 521439867$.

The map $\phi$ defines a folklore bijection between noncrossing partitions of size $n$ with $l+1$ blocks and $231$-avoiding permutations with $l$ descents. This recovers known numerology about the $(0,2)$-case.

\begin{remark}
    More generally, standard segmented permutations with zero $\sminv$ can be characterized as $231$-avoiding permutations where letters of a block are smaller all than letters of the blocks to its right. Of these, there are $\binom{2n+1}{n}$ as is easily shown, and thus we recover the total dimension of the $(0,2)$-coinvariant ring.
\end{remark}

\bibliographystyle{plain}
\bibliography{references}

\end{document}